\documentclass{amsart}
\usepackage{amsmath,amssymb,amsxtra,amsthm,amscd}
\usepackage{lmodern}
\usepackage[all,cmtip]{xy}
\usepackage{rotating}
\usepackage{microtype}
\usepackage{todonotes}

\usepackage{adjustbox}
\usepackage{color}

\newif\ifShowLabels
\ShowLabelstrue
\newcommand{\TeXref}[1]{
\marginpar{\scriptsize \texttt{#1}}}

\ShowLabelsfalse

\DeclareMathOperator{\EGamma}{\boldsymbol{E}\boldsymbol{\Gamma}}

\DeclareMathOperator{\fil}{fil}

\DeclareMathOperator{\Free}{\mathbf{Free}}

\DeclareMathOperator{\Fun}{Fun}

\DeclareMathOperator{\Hom}{Hom}

\DeclareMathOperator{\id}{id}

\DeclareMathOperator{\im}{im}

\DeclareMathOperator{\K}{\mathit{K}}
         \newcommand{\Knc}{\K^{-\infty}}

         \newcommand{\KncX}{\K^{-\infty}_X}

\DeclareMathOperator{\point}{pt}

\DeclareMathOperator{\Sh}{Sh}

\DeclareMathOperator*{\one}{1}
\newcommand{\onehatplace}[1]
{ \one^{\substack{#1 \\ \frown}} }

\DeclareMathOperator*{\bones}{\times}
\newcommand{\undertimes}[1]
{ \bones_{#1} }

\DeclareMathOperator*{\bowl}{\cup}
\newcommand{\undercup}[1]
{ \bowl_{#1} }

\DeclareMathOperator*{\arch}{\cap}
\newcommand{\undercap}[1]
{ \arch_{#1} }

\newcommand{\pull}
{\!\!\! -\!\!\! -\!\!\! -\!\!\!}

\DeclareMathOperator*{\holimprep}{holim}                       
\newcommand{\holim}[1]%
{\displaystyle\holimprep_{\substack{\leftarrow \pull - \\ #1}} \, }

\DeclareMathOperator*{\hocolimprep}{hocolim}                   
\newcommand{\hocolim}[1]%
{\displaystyle\hocolimprep_{\substack{- \pull \rightarrow \\ #1}} \, }

\DeclareMathOperator*{\plainlim}{lim}                           
\newcommand{\contralim}[1]%
{\displaystyle\plainlim_{\substack{\leftarrow \pull - \\ #1}} \, }

\DeclareMathOperator*{\plaincolim}{colim}                       
\newcommand{\colim}[1]%
{\displaystyle\plaincolim_{\substack{- \pull \rightarrow \\ #1}} \, }

\DeclareMathOperator*{\laxlimplain}{laxlim}                     
\newcommand{\laxlim}[1]%
{\displaystyle\laxlimplain_{\substack{\leftarrow \pull - \\ #1}} \, }

\providecommand{\bysame}{\makebox[3em]{\hrulefill}\thinspace}

\theoremstyle{plain}
\newtheorem{Thm}{Theorem}[section]

\newtheorem{Cor}[Thm]{Corollary}

\newtheorem{Lem}[Thm]{Lemma}
\newtheorem{Prop}[Thm]{Proposition}

\theoremstyle{definition}
\newtheorem{Def}[Thm]{Definition}

\newtheorem{Ex}[Thm]{Example}
\newtheorem{Exs}[Thm]{Examples}

\newtheorem{Rem}[Thm]{Remark}

\newtheorem{Que}[Thm]{Question}
\newtheorem{Conj}[Thm]{Conjecture}

\theoremstyle{remark}
\newtheorem{Not}[Thm]{Notation}

\newtheoremstyle{freestylethm}{6pt}{6pt}{\itshape}{}%
                {\bfseries}{}{.5em}{\thmnote{#3}}
\theoremstyle{freestylethm}

\newcommand{\SecRef}[2]{\section{#1}\label{S:#2}%
\ifShowLabels \TeXref{{S:#2}} \fi}

\newcommand{\refS}[1]{\textup{\ref{S:#1}}}

\newcommand{\refT}[1]{\textup{\ref{T:#1}}}
\newcommand{\refL}[1]{\textup{\ref{L:#1}}}
\newcommand{\refD}[1]{\textup{\ref{D:#1}}}
\newcommand{\refC}[1]{\textup{\ref{C:#1}}}

\newcommand{\refP}[1]{\textup{\ref{P:#1}}}
\newcommand{\refR}[1]{\textup{\ref{R:#1}}}
\newcommand{\refN}[1]{\textup{\ref{N:#1}}}

\newenvironment{ThmRef}[1]%
{ \begin{Thm} \label{T:#1}
\ifShowLabels \TeXref{T:#1} \fi }%
{ \end{Thm} }
\newenvironment{DefRef}[1]%
{ \begin{Def} \label{D:#1}
\ifShowLabels \TeXref{D:#1} \fi }%
{ \end{Def} }
\newenvironment{LemRef}[1]%
{ \begin{Lem} \label{L:#1}
\ifShowLabels \TeXref{L:#1} \fi }%
{ \end{Lem} }
\newenvironment{CorRef}[1]%
{ \begin{Cor} \label{C:#1}
\ifShowLabels \TeXref{C:#1} \fi }%
{ \end{Cor} }
\newenvironment{RemRef}[1]%
{ \begin{Rem} \label{R:#1}
\ifShowLabels \TeXref{R:#1} \fi }%
{ \end{Rem} }
\newenvironment{PropRef}[1]%
{ \begin{Prop} \label{P:#1}
\ifShowLabels \TeXref{P:#1} \fi }%
{ \end{Prop} }
{ \begin{Ex} \label{E:#1}
\ifShowLabels \TeXref{E:#1} \fi  }%
{ \end{Ex} }
\newenvironment{NotRef}[1]%
{ \begin{Not} \label{N:#1}
\ifShowLabels \TeXref{N:#1} \fi }%
{ \end{Not} }
{ \begin{Que} \label{Q:#1}
\ifShowLabels \TeXref{Q:#1} \fi }%
{ \end{Que} }
{ \begin{Conj} \label{Conj:#1}
\ifShowLabels \TeXref{Conj:#1} \fi }%
{ \end{Conj} }

\newenvironment{ThmRefName}[2]%
{ \begin{Thm} [#2]\label{T:#1}
\ifShowLabels \TeXref{T:#1} \fi }%
{ \end{Thm} }
\newenvironment{DefRefName}[2]%
{ \begin{Def} [#2]\label{D:#1}
\ifShowLabels \TeXref{D:#1} \fi }%
{ \end{Def} }
{ \begin{Lem} [#2]\label{L:#1}
\ifShowLabels \TeXref{L:#1} \fi }%
{ \end{Lem} }
{ \begin{Cor} [#2]\label{C:#1}
\ifShowLabels \TeXref{C:#1} \fi }%
{ \end{Cor} }
{ \begin{Rem} [#2]\label{R:#1}
\ifShowLabels \TeXref{R:#1} \fi }%
{ \end{Rem} }
{ \begin{Prop} [#2]\label{P:#1}
\ifShowLabels \TeXref{P:#1} \fi }%
{ \end{Prop} }
{ \begin{Ex} [#2]\label{E:#1}
\ifShowLabels \TeXref{E:#1} \fi }%
{ \end{Ex} }

\let\oldtocsection=\tocsection
\let\oldtocsubsection=\tocsubsection
\renewcommand{\tocsection}[2]{\hspace{0em}\oldtocsection{#1}{#2}}
\renewcommand{\tocsubsection}[2]{\hspace{2em}\oldtocsubsection{#1}{#2}}

\begin{document}

\title{Asymptotic transfer maps in parametrized \textit{K}-theory}
\author[Gunnar Carlsson]{Gunnar Carlsson}
\address{Department of Mathematics\\ Stanford University\\ Stanford\\ CA 94305}
\email{gunnar@math.stanford.edu}
\author[Boris Goldfarb]{Boris Goldfarb}
\address{Department of Mathematics and Statistics\\ SUNY\\ Albany\\ NY 12222}
\email{goldfarb@math.albany.edu}
\date{\today}

\begin{abstract}
We define asymptotic transfers in bounded $K$-theory together with a context where this can be done in great generality.  Controlled algebra plays a central role in many advances in geometric topology, including recent work on Novikov, Borel, and Farrell-Jones conjectures.  One of the features that appears in various manifestations throughout the subject, starting with the original work of Farrell and Jones, is an asymptotic transfer whose meaning and construction depend on the geometric circumstances.  We first develop a general framework that allows us to construct a version of asymptotic transfer maps for any finite aspherical complex.  This framework is the equivariant parametrized $K$-theory with fibred control.  We also include several fibrewise excision theorems for its computation and a discussion of where the standard tools break down and which tools replace them.  
\end{abstract}

\maketitle

\tableofcontents

\SecRef{Introduction}{INTRO}

Controlled algebra of geometric modules uses geometric control conditions well-suited for generating $K$-theory spectra out of the associated additive controlled categories of free modules.  Specifically, the bounded controlled algebra was developed by E.K. Pedersen and C. Weibel in \cite{ePcW:85,ePcW:89}.  It has been used extensively for delooping $K$-theory of rings and other algebraic objects and in the study of the Novikov conjecture \cite{gC:95,gCbG:04,gCeP:93} about the assembly map for algebraic $K$-theory and $L$-theory of certain group rings. 

There are two innovations in this paper: the setting of parametrized $K$-theory with fibred control which generalizes bounded $K$-theory and a general asymptotic transfer that can be defined in this new theory.
   
We start by explaining the work of Pedersen/Weibel as an algebraic theory of free modules, possibly infinitely generated, parametrized by a metric space.  The input is a proper metric space $X$ (i.e. a metric space in which every closed bounded subset is compact) and a ring $R$.  Given this data, consider triples $(M, B, \varphi)$, where $M$ is a free left $R$-module with basis $B$, and where $\varphi \colon B \rightarrow X$ is a reference function with the property that the inverse image of any bounded subset is finite.  A morphism from $(M,B,\varphi)$ to $(M^{\prime}, B^{\prime} , \varphi^{\prime})$  is an $R$-module homomorphism $f$ from $M$ to $M^{\prime}$  which has the property that there exists a bound $b \ge 0$ so that for any basis element $\beta \in B$, $f(\beta)$ is in the span of basis elements $\beta' \in B^{\prime}$ for which $d (\varphi (\beta), \varphi' (\beta ^{\prime})) \leq b$.  This is the simplest version of a {\em control condition} that can be imposed on homomorphisms to construct various categories of modules, and it leads to the category of \textit{geometric modules} $\mathcal{C} (X, R)$.  This is an additive category to which one can apply the usual algebraic $K$-theory construction.   

A formulation of this theory that was observed already in \cite{ePcW:85} allows one to use objects from an arbitrary additive category $\mathcal{A}$ as ``coefficients'' in place of the finitely generated free $R$-modules.  We will spell out the details of this further in the paper.  So one obtains a new additive category $\mathcal{C} (X, \mathcal{A})$. This observation leads to the possibility of iterating the geometric control construction: if there are two proper metric spaces $X$ and $Y$, then we have an additive category $\mathcal{C} (X, \mathcal{C} (Y, R))$.  This is precisely a case of \textit{fibred control} we want to study in this paper. 

To give an idea for what kind of control is implied by this construction, let's parse the implications for objects viewed as $R$-modules.  In this setting an object is a free module $M$ viewed as parametrized over the product $X \times Y$, so we have a pair $(M,B)$ with a reference function $\varphi \colon B \rightarrow X \times Y$.  The control condition on homomorphisms $f$ from $(M_1, B_1, \varphi _1)$ to $(M_2, B_2, \varphi _2)$ amounts to the existence of a number $b_X \ge 0$ and a function $c \colon X \to [0, \infty)$ such that for $\beta \in B_1$, we have that $f(\beta)$ is a linear combination of basis elements $\beta ^{\prime} \in B_2$ so that $d(\pi _X(\varphi _1(\beta)), \pi _X(\varphi _2 (\beta ^{\prime}))) \leq b_X$ and 
 $d(\pi_Y (\varphi _1(\beta)), \pi_Y (\varphi _2 (\beta ^{\prime}))) \leq c (\pi _X(\varphi _1(\beta)))$.  
Contrast this with the control condition in $\mathcal{C} (X \times Y, R)$ which is exactly as above except $c$ can be chosen to be a constant function.
This should suggest that the condition of fibred control, with $X$ regarded as a base and $Y$ regarded as a fibre, relaxes the usual control condition over $X \times Y$ by letting $c$ be a quantity varying with $x \in X$. 

It is shown in \cite{gC:95} that $\mathcal{C} (X, R)$ enjoys a certain excision property, which permits among other things a comparison with Borel-Moore homology spectra with coefficients in the $K$-theory spectrum $K(R)$. 
This property becomes crucial in the equivariant homotopy theoretic approach to the Novikov conjecture.
In this paper we prove a number of excision results in categories with fibred control and some related categories of fixed-point objects with respect to naturally occurring actions on the metric spaces.
These theorems are required for further work on assembly maps in $K$-theory.

We briefly address the methods used and the interesting phenomena that come up.  The go-to technique for proving localization and excision theorems in controlled algebra is the use of Karoubi filtrations in additive categories.  This technique is indeed used also in our proofs of general non-equivariant fibrewise excision results (section \refS{KTFCMNB}) and in some specific equivariant situations (section \refS{JKASEQT}).  However, we want point out one curious situation that will interest the experts. The most useful novelty of the fibred control is that it allows introduction of constraints on features, including actions, that vary across fibres.  When this is done to actions, it becomes unnatural to insist that the action is by maps that are necessarily isometries on the nose or, if the action is by more general coarse equivalences, that the same numerical constraint is satisfied across all fibres.  In this situation we can state desired excision statements about fixed point object categories for multiple actions but discover that they are inaccessible to the Karoubi filtration technique.  We give an example in section \refS{MExD} where we are able to pinpoint the deficiency that makes the Karoubi filtrations unavailable. 

Before we move on to the asymptotic transfer, we want to point out that the generalization of controlled algebra in this paper is useful for a variety of other applications.  
It is most immediately the natural controlled theory to consider for bundles on non-compact manifolds where the control in the fibre direction can vary with the fibre.
This holds more generally for stacks or the coarse analogues of bundles that have already appeared in the literature \cite{kW:10}. 
Another manifestation of fibred control is found in the active current research area seeking sufficient geometric conditions for verifying the coarse Baum-Connes conjecture.  The original condition in \cite{gY:98,gY:00}  was existence of a coarse embedding of the group into a Hilbert space. The embedding can be viewed as a way to gain geometric control for performing stable excision over the group or the space with a proper cocompact action by the group.   That condition has been relaxed by several authors to \textit{fibred coarse embeddings} \cite{sA:16,xCqWxW:13,xCqWgY:13,mF:14,mMhS:18}.  The results of this paper allow to address the Borel conjecture in $K$-theory for groups admitting a fibred coarse embedding into Hilbert space.

A transfer map is a very useful tool in algebraic topology. For a continuous map between spaces $f \colon X \to Y$ and a covariant functor $F$ on spaces to some algebraic category, a transfer $T (f)$ is a morphism $F(Y) \to F(X)$, usually with several properties that depend on the context.  This is a ``wrong way'' map as opposed to the usual induced maps $F(f) \colon F(X) \to F(Y)$.  A survey of the general homotopy theoretic construction for projections in Hurewicz fibrations with homotopy finite fibers and numerous points of views and applications is given by its creators in Becker/Gottlieb \cite{jBdG:99}.  Infinite transfers with no finiteness assumption on the fibers appear naturally in the area of topology related to modern proofs of the Novikov and Farrell-Jones conjectures.  Here the functor $F$ is usually either locally finite homology or controlled $K$-theory.  We refer to \cite{gC:93} for a general introduction to the ideas and to \cite{gC:95} for the details required for the proofs of the Novikov conjecture about injectivity of assembly maps in $K$-theory.  In this paper, we define a parametrized version of the infinite transfer and establish some of its properties needed for the study of surjectivity of those maps.  

Asymptotic transfers are a crucial feature in modern proofs of cases of the Borel rigidity conjecture. This is most notable in the pioneering work of Farrell and Jones \cite{fFlJ:86}. The basic  idea in controlled topology is that improvement in control leads to trivialization of an appropriate h-cobordism or a surgery problem. In the case of closed hyperbolic manifolds, Farrell and Jones were able to define a special asymptotic transfer to the unit sphere bundle of the manifold and use the geodesic flow on the bundle to control the lifted homotopies and eventually trivialize the obstructions to geometric problems.
Further work on the Farrell-Jones conjecture relied in similar ways on an asymptotic transfer to a setting where a given algebraic representation of the geometric problem can be manipulated to have progressively better control. This can be done in cases where subtle consequences of immanent manifestations of non-positive curvature can be exploited for controlling geometric problems.  We refer to instances of this in section 6 of \cite{aBwLhR:08a} or section 7 of \cite{aBwL:12}.  In contrast to these cases, we construct in this paper a $K$-theory transfer that uses no geometric constraints.  The target of this transfer is the parametrized $K$-theory where, in analogous fashion, a particular trivialization exists as indicated at various places, for example in Theorem \refT{GpDef}.  For the application to aspherical manifolds, the fibre for the parametrized $K$-theory is the universal cover of the normal disk bundle to the embedding of the manifold in a Euclidean space.  It is a metric space with the metric constructed in section \refS{PREPSEC}.  Trivialization in parametrized $K$-theory, after certain algebraic constraints are imposed on the coefficient ring, is an outcome of a general equivariant excision theorem proven elsewhere \cite{gCbG:20}.  We address the relation between that $G$-theoretic excision theorem and the excision theorems from section \refS{JKASEQT} in Remark \refR{mvhfbd}.

The paper starts with a paced introduction to bounded $K$-theory and Karoubi filtration methods used to compute it.  After that, we introduce fibred control and the equivariant fibred $K$-theory in sections \refS{KTFCMNB} and \refS{JKASEQT}.  The material from section \refS{BETRGCV} is used to motivate and construct the proofs of similar controlled excision theorems in the parametrized setting.  We also use it to carefully describe in section \refS{MExD} how the standard techniques fail in proving a natural desired statement in the parametrized theory.  This motivates the development of fibred $G$-theory in \cite{gCbG:19,gCbG:20} where the equivariant fibred excision theorems do hold. Finally, the last two sections \refS{PREPSEC} and \refS{CPT} are about the asymptotic transfer.

\SecRef{Pedersen-Weibel categories, Karoubi filtrations}{BETRGCV}

Bounded $K$-theory introduced in Pedersen \cite{eP:84} and Pedersen--Weibel \cite{ePcW:85} associates a nonconnective spectrum $\Knc (M,R)$ to a proper metric space $M$ (a metric space where closed bounded subsets are compact) and an associative ring $R$ with unity.  We are going to start with a careful, at times revisionist, review of the well-known features of this theory that will need generalization. 

The metric spaces in this subject are often understood in the generalized sense as follows.  

\begin{DefRef}{GenMS}
  A \textit{generalized metric space} is a set $X$ and a function
$d \colon X \times X \to [0,\infty) \cup \{\infty\}$
which is reflexive, symmetric, and satisfies the triangle inequality
in the usual way.
Classical metric spaces are the generalized metric spaces with distance function $d$ assuming only finite values.	
We will use the term \textit{metric space} to mean a generalized metric space.

We call a subset $M'$ of a metric space $M$ a \textit{metric subset} if all values of $d \vert M' \times M'$ are finite.  The maximal metric subsets are called \textit{metric components}.
A metric space is \textit{proper} if it is a countable disjoint union of
metric components $M_i$,
and all closed metric balls in $M$ are compact.
The metric topology on a metric space is
defined as usual.  
\end{DefRef}

In this section, all metric spaces and facts about them may be viewed in this generalized context.  We will restrict to spaces with a single metric component starting in the next section.

\begin{DefRef}{PWcats}
$\mathcal{C} (M,R)$ is the additive category of \textit{geometric $R$-modules} whose 
objects are functions $F \colon M \to \Free_{fg} (R)$ which are locally finite assignments of free finitely generated $R$-modules $F_m$ to points $m$ of $M$.
The local finiteness condition requires precisely that for any bounded subset $S \subset M$ the restriction of $F$ to $S$ has finitely many nonzero values.
Let $d$ be the distance function in $M$.  The morphisms in $\mathcal{C} (M,R)$ are the $R$-linear homomorphisms
\[
\phi \colon \bigoplus_{m \in M} F_m \longrightarrow \bigoplus_{n \in M} G_n
\]
with the property
that the components $F_m \to G_n$ are zero for $d(m,n) > b$
for some fixed real number $b = b (\phi) \ge 0$.
The associated $K$-theory spectrum is denoted by $K (M,R)$, or $K (M)$ when the choice of ring $R$ is implicit,
and is called the \textit{bounded K-theory} of $M$.
\end{DefRef}

\begin{NotRef}{ttyh}
For a subset $S \subset M$ and a real number $r \ge 0$,
$S[r]$ will stand for the metric $r$-enlargement $\{ m \in M \mid d (m,S) \le r \}$.
In this notation, the metric ball of radius $r$ centered at $x$ is $\{ m \} [r]$ or simply $m [r]$.
\end{NotRef}

Now for every object $F$ and subset $S$ there is a free $R$-module $F(S) = \bigoplus_{m \in S} F_m$.  The condition that $\phi$ is controlled as above is equivalent to existence of a number $b \ge 0$ so that $\phi F(S) \subset F(S[b])$ for all choices of $S$.

\begin{PropRef}{KLOWAS}
	The description of $\mathcal{C} (M,R)$ in the introduction defines a category additively equivalent to the one in Definition \refD{PWcats}, establishing a dictionary between terminology in various papers in the literature.
\end{PropRef}

\begin{proof}
	Given a triple $(M, B, \varphi)$ described in the introduction, define $M_x$ as the $R$-submodule freely generated by $\varphi^{-1} (x)$.  It is clear that this gives an additive functor in one direction.  The inverse functor is constructed by selecting a finite basis in each $F_x$ and defining $B$ to be the union of these bases.  The map $\varphi$ sends a basis element $b$ to $x$ if $b \in F_x$.
\end{proof}

Inclusions of metric spaces induce additive functors between the corresponding bounded $K$-theory spectra.
The main result of Pedersen--Weibel \cite{ePcW:85} is a delooping theorem which can be stated as follows.

\begin{ThmRefName}{Deloop}{Nonconnective delooping of bounded \textit{K}-theory}
Given a proper metric space $M$ and the standard Euclidean metric on the real line $\mathbb{R}$, the natural inclusion $M \to M \times \mathbb{R}$ induces isomorphisms $K_n (M) \simeq K_{n-1} (M \times \mathbb{R})$
for all integers $n > 1$.
If one defines the spectrum
\[
\Knc (M,R) \ = \ \hocolim{k} \Omega^k K (M \times \mathbb{R}^k),
\]
then the stable homotopy groups of  $\Knc (R) = \Knc (\point,R)$ coincide with the algebraic $K$-groups of $R$ in positive dimensions and with the Bass negative $K$-theory of $R$ in negative dimensions.
\end{ThmRefName}

When we develop the equivariant theory, we will want to consider group actions by maps that are more general than isometries.
Let $X$ and $Y$ be proper metric spaces with metric functions $d_X$ and $d_Y$.
This means, in particular, that closed bounded subsets of $X$ and $Y$ are compact.

\begin{DefRef}{Bornol2}
A function $f \colon X \to Y$
between proper metric spaces is \textit{uniformly expansive} if there is a real positive function $l$
such that for pairs of points $x_1$ and $x_2$ the inequality 
$d_X (x_1, x_2) \le r$ implies the inequality $d_Y (f(x_1), f(x_2)) \le l(r)$. 
This is the same concept as \textit{bornologous} maps in Roe \cite[Definition 1.8]{jR:03}.  It should be emphasized that the function is not assumed to be necessarily continuous.  
The function $f$ is \textit{proper} if $f^{-1} (S)$ is a bounded subset of $X$ for
each bounded subset $S$ of $Y$.
We say $f$ is a \textit{coarse map} if it is proper and uniformly expansive.
\end{DefRef}

\begin{ThmRef}{JARI}
Coarse maps between proper metric spaces induce additive functors between bounded categories.  
\end{ThmRef}

We include a proof of this basic fact that will be generalized in parametrized setting.

\begin{proof}
Let $f \colon X \to Y$ be coarse.
The additive functor $f_{\ast} \colon \mathcal{C} (X,R) \to \mathcal{C} (Y,R)$.
is induced on objects by the assignment
\[
(f_{\ast} F)_y = \bigoplus_{x \in f^{-1} (y)} F_x.
\]
Since $f$ is proper, $f^{-1} (y)$ is a bounded set for all $y$ in $Y$.
So the direct sum in the formula is finite, and $(f_{\ast} F)_y$ is a finitely generated free $R$-module.
If $S \subset Y$ is a bounded subset then $f^{-1} (S)$ is bounded.
There are finitely many $F_z \ne 0$  for $z \in f^{-1} (S)$ and therefore finitely many $(f_{\ast} F)_y \ne 0$ for $y \in S$.
This shows $f_{\ast} F$ is locally finite.

Notice that
\[
f_{\ast} F = \bigoplus_{y \in Y} (f_{\ast} F)_y = \bigoplus_{y \in Y} \bigoplus_{z \in f^{-1} (y)} F_z = F.
\]

Suppose we are given a morphism $\phi \colon F \to G$ in $\mathcal{C} (X,R)$.
Interpreting $f_{\ast} F$ and $f_{\ast} G$ as the same $R$-modules, as in the formula above, we define $f_{\ast} \phi \colon f_{\ast} F \to f_{\ast} G$ equal to $\phi$.
We must check that $f_{\ast} \phi$ is bounded.
Suppose $\phi$ is bounded by $D$, and $f$ is $l$-coarse.  We claim that $f_{\ast} \phi$ is bounded by $l (D)$.  Indeed, if $d_{Y} (y,y') > l (D)$ then $d_X (x,x') > D$ for all $x$, $x' \in X$ such that $f(x) = y$ and $f(x') = y'$.  So all components $\phi_{x,x'} = 0$, therefore all components $(f_{\ast} \phi)_{y,y'} = 0$.
\end{proof}

\begin{CorRef}{HMLOWQ}
$\Knc$ is a covariant functor from the category of proper metric spaces and coarse maps to the category of spectra.
\end{CorRef}

The map $f$ is a \textit{coarse equivalence} if there is a coarse map $g \colon Y \to X$ such
that $f \circ g$ and $g \circ f$ are bounded maps.

\begin{Exs} \label{POIU}
The isometric embedding of a metric subspace is a coarse map.
An isometry, which is a bijective isometric map, is a coarse equivalence.
An isometric embedding onto a subspace that has the property that its bounded enlargement is the whole target metric space is also a coarse equivalence.

Any \textit{bounded} function $f \colon X \to X$, that is a function with
$d_X (x, f(x)) \le D$ for all $x \in X$ and a fixed $D \ge 0$, is a coarse equivalence.
More generally, a quasi-isometry $f \colon X \to Y$ onto a subset $U \subset Y$ such that for some number $s \ge 0$ we have $U[s]=Y$ is a coarse equivalence.
\end{Exs}

The following definition makes precise a useful class of metrics one has on a finitely generated group.

\begin{DefRef}{WORD}
The \textit{word-length metric} $d = d_{\Omega}$ on a group
$\Gamma$ with a fixed finite generating set $\Omega$ closed under taking inverses is the length metric
induced from the condition that $d (\gamma, \gamma \omega) =1$, whenever $\gamma \in \Gamma$ and $\omega \in \Omega$.  In other words, $d (\alpha, \beta)$ is the minimal length $t$ of sequences $\alpha = \gamma_0, \gamma_1, \ldots  \gamma_t = \beta$ in $\Gamma$ where each consecutive pair of elements differs by right multiplication by an $\omega$ from $\Omega$.
This metric makes $\Gamma$ a proper metric space with a free action by $\Gamma$ via left
multiplication.
\end{DefRef}

If one considers a different choice of a finite generating set $\Omega'$, it is well-known that the identity map on the group with the two metrics $d_{\Omega}$ and $d_{\Omega'}$ is a quasi-isometry, cf. \cite[Proposition 1.15]{jR:03}.  We see from the combination of Corollary \refC{HMLOWQ} and the examples above that the bounded $K$-theory of $\Gamma$ is independent from the choice of a finite generating set, up to an equivalence.

Our special interest in this paper is in actions by bounded coarse equivalences.  

\begin{DefRef}{leftbdd}
A left action of $\Gamma$ on a metric space $X$ is \textit{bounded} if for each element $\gamma \in \Gamma$ there is a number $B_{\gamma} \ge 0$ such that $d(x, \gamma x) \le B_{\gamma}$ for all $x \in X$.
\end{DefRef}

Such actions do appear naturally but only in special circumstances.  If $\Gamma$ is a group acting on itself via left multiplication, then this action is a bounded action by isometries when the group is a finitely generated abelian group.  However, the left multiplication action of the integral lattice in the three-dimensional Heisenberg group is already not bounded because of the well-known warping of the metric.

In section \refS{PREPSEC} we will define a useful conversion of any action by isometries to a bounded action which is sufficient for our and many other purposes.  

\medskip

We next review an excision theorem that makes bounded $K$-theory computable in special but crucial geometric situations.  This review includes some details that we will use to prove some excision results for parametrized $K$-theory in the next section.

Suppose $U$ is a subset of $M$.
Let $\mathcal{C} (M,R)_{<U}$ denote the full subcategory of $\mathcal{C} (M,R)$ on the objects $F$ with $F_m = 0$ for all points $m \in M$ with $d (m, U) > D$ for some fixed number $D > 0$ specific to $F$.
This is an additive subcategory of $\mathcal{C} (M,R)$ with the associated $K$-theory spectrum $\Knc (M,R)_{<U}$.

Similarly, if $U$ and $V$ are a pair of subsets of $M$, then there is the full additive subcategory
$\mathcal{C} (M,R)_{<U,V}$ of $F$ with $F_m = 0$ for all $m$ with $d (m, U) > D_1$ and $d (m, V) > D_2$ for some numbers $D_1, D_2 \ge  0$.
It is easy to see that $\mathcal{C} (U,R)$ is a subcategory which is in fact equivalent to $\mathcal{C} (M,R)_{<U}$.

\medskip

The following theorem is the basic computational device in bounded $K$-theory.

\begin{ThmRefName}{BddExc}{Bounded Excision, Theorem IV.1 of \cite{gC:95}}
Given a proper metric space $M$ and a pair of subsets $U$, $V$ of $M$ which form a cover of $M$, there is a homotopy pushout diagram 
\[
\xymatrix{
 \Knc (M)_{<U, V} \ar[r] \ar[d]
&\Knc (U) \ar[d] \\
 \Knc (V) \ar[r]
&\Knc (M) }
\]
\end{ThmRefName}

It is possible to restate the Bounded Excision Theorem in a more intrinsic form, after restricting to a special class of coverings.

\begin{DefRef}{CATHJK}
A pair of subsets $S$, $T$ of a proper metric space $M$ is called \textit{coarsely antithetic} if $S$ and $T$ are proper metric subspaces with the subspace metric and for each pair of numbers $d_S$, $d_T \ge  0$ there is a number $d' \ge  0$ so that
\[
S[d_S] \cap T[d_T] \subset (S \cap T)[d'].
\]
\end{DefRef}

Examples of coarsely antithetic pairs
include complementary closed half-spaces in a Euclidean space, as well as
any two non-vacuously intersecting closed subsets of a simplicial tree.
In the latter example, a tree is viewed as a geodesic metric space where the metric is induced from the local condition that each edge is isometric to a closed real interval.

\begin{CorRef}{BDDEXCII}
If $U$ and $V$ is a coarsely antithetic pair of subsets of $M$ which form a cover of $M$, then the commutative square
\[
\xymatrix{
 \Knc (U \cap V) \ar[r] \ar[d]
&\Knc (U) \ar[d] \\
 \Knc (V) \ar[r]
&\Knc (M)
}
\]
is a homotopy pushout.
\end{CorRef}

We want to outline the proof of Theorem \refT{BddExc} in specific terms that will be used later.  Another reason for recording rather standard details is to refer to them later when we study a failure of these methods in section \refS{MExD}.

The notion of Karoubi filtrations in additive categories is central to the proof of this theorem as developed by Cardenas/Pedersen \cite{mCeP:97}.

\begin{DefRef}{Karoubi}
An additive category $\mathcal{C}$ is \textit{Karoubi filtered} by a full subcategory $\mathcal{A}$
if every object $C$ of $\mathcal{C}$ has a family of decompositions $\{ C = E_{\alpha} \oplus D_{\alpha} \}$ with $E_{\alpha} \in \mathcal{A}$ and $D_{\alpha} \in \mathcal{C}$, called a \textit{Karoubi filtration} of $C$, satisfying the following properties.
\begin{itemize}
\item For each object $C$ of $\mathcal{C}$, there is a partial order on Karoubi decompositions such that $E_{\alpha} \oplus D_{\alpha} \le E_{\beta} \oplus D_{\beta}$ whenever $D_{\beta} \subset D_{\alpha}$ and $E_{\alpha} \subset E_{\beta}$.
\item Every morphism $A \to C$ with $A \in \mathcal{A}$ and $C \in \mathcal{C}$ factors as $A \to  E_{\alpha} \to E_{\alpha} \oplus D_{\alpha} = C$ for some value of $\alpha$.
\item Every morphism $C \to A$  with $C \in \mathcal{C}$ and $A \in \mathcal{A}$ factors as $C = E_{\alpha} \oplus D_{\alpha} \to  E_{\alpha} \to A$ for some value of $\alpha$.
\item For each pair of objects $C$ and $C'$ with the corresponding filtrations $\{ E_{\alpha} \oplus D_{\alpha} \}$ and $\{ E'_{\alpha} \oplus D'_{\alpha} \}$, the filtration of $C \oplus C'$ is the family $\{ C \oplus C' = ( E_{\alpha} \oplus E'_{\alpha} ) \oplus ( D_{\alpha} \oplus D'_{\alpha} ) \}$.
\end{itemize}
A morphism $f \colon C \to D$ in $\mathcal{C}$ is $\mathcal{A}$-\textit{zero} if $f$ factors through an object of $\mathcal{A}$.
One defines the \textit{Karoubi quotient} $\mathcal{C}/\mathcal{A}$ as the additive category with the same objects as $\mathcal{C}$ and morphism sets $\Hom_{\mathcal{C}/\mathcal{A}} (C,D) = \Hom (C,D)/ \{ \mathcal{A}\mathrm{-zero \ morphisms} \}$.
\end{DefRef}

The following is the main theorem of Cardenas--Pedersen \cite[Theorem 7.1]{mCeP:97}.

\begin{ThmRefName}{CPF}{Fibration Theorem}
Suppose $\mathcal{C}$ is an $\mathcal{A}$-filtered category, then there is a homotopy fibration
\[
K ( \mathcal{A}^{\wedge K} ) \longrightarrow
K ( \mathcal{C} ) \longrightarrow K ( \mathcal{C}/\mathcal{A} ).
\]
Here $\mathcal{A}^{\wedge K}$ is a certain subcategory of the idempotent completion of $\mathcal{A}$ with the same positive $K$-theory as $\mathcal{A}$.
\end{ThmRefName}

The next statement is a consequence of the Fibration Theorem obtained in the last paragraph of \cite{mCeP:97}.

\begin{CorRef}{CPFrev}
Suppose $\mathcal{C}$ is an $\mathcal{A}$-filtered category, then there is a homotopy fibration
\[
\Knc ( \mathcal{A} ) \longrightarrow
\Knc ( \mathcal{C} ) \longrightarrow \Knc ( \mathcal{C}/\mathcal{A} ).
\]
\end{CorRef}

Theorem \refT{BddExc} follows from this Corollary by the following device.

\begin{proof}[Sketch of the proof of Theorem \refT{BddExc}]
The first crucial observation is that $\mathcal{C} = \mathcal{C}(M)$ is $\mathcal{A} = \mathcal{C} (M)_{< U}$-filtered.
For future references, let us spell out what is involved.

The additive structure in $\mathcal{C}$ is given by $(F \oplus G)_m = F_m \oplus G_m$.  So, in particular, $(F \oplus G)(S) = F(S) \oplus G(S)$ for all subsets $S$.
Given an object $F$ of $\mathcal{C}$, the subobjects $F (U[k])$, for $k \ge 0$, are free direct summands of $F(M)$ which can be given the structure of a geometric module over $M$ in the obvious way.  Now the decompositions 
\[
F = F (U[k]) \oplus F (M \setminus U[k]) 
\]
is the family we need for a Karoubi filtration.  Suppose, for illustration, we have $f \colon A \to F$ bounded by $b \ge 0$.  Then 
$A(M) = A(U[r])$ for some $r \ge 0$ so $f (A) \subset F(U[r + b])$.  So, indeed, $f$ factors through this direct summand.

Since $\mathcal{C}(M)$ is $\mathcal{C} (M)_{< V}$-filtered, there is the additive Karoubi quotient which we denote by $\mathcal{C} (M,V)$, with the nonconnective $K$-theory $\Knc (M, V)$.  

For simplicity, let us assume that $U$ and $V$ form a coarsely antithetic pair, then $\mathcal{C} (M)_{< U}$ is similarly $\mathcal{C} (M)_{< U,V}$-filtered, with the Karoubi quotient $\mathcal{C} (U, U \cap V)$.
Corollary \refC{CPFrev} gives two homotopy fibrations that form a commutative diagram

\[
\begin{CD}
\Knc (M)_{< U,V} @>>> \Knc (M)_{< U} @>>> \Knc (U, U \cap V) \\
@VVV @VVV @VV{\simeq}V \\
\Knc (M)_{< V} @>>> \Knc (M) @>>> \Knc (M, V)
\end{CD}
\]
\\ \indent
The second crucial observation is that the exact inclusion of additive categories 
\begin{equation}
	I \colon \mathcal{C} (U, U \cap V) \longrightarrow  \mathcal{C} (M,V) \tag{\dag} \label{ISOM}
\end{equation}
is an equivalence. Indeed, the restriction which sends an object $F$ to the direct summand $F(U)$ gives a functor $J \colon \mathcal{C} (M,V) \to \mathcal{C} (U, U \cap V)$. The composition $J \circ I$ is the identity, while $I \circ J$ possesses a natural transformation to the identity where each morphism is the projection $F \to F(U)$. Since the complement of $U$ is contained in $V$, this projection is an isomorphism in $\mathcal{C} (M,V)$.  This equivalence of categories induces a weak equivalence in $K$-theory and proves the theorem. 
\end{proof}

\SecRef{Parametrized $K$-theory with fibred control}{KTFCMNB}

When $M$ is the product of two proper metric spaces $X \times Y$ with the metric
\[
  d_{\textrm{max}} ((x_1,y_1),(x_2,y_2)) = \max \{ d_X (x_1,x_2),
  d_Y (y_1,y_2) \},
\]
one has the bounded category of geometric $R$-modules $\mathcal{C} (X \times Y, R)$ of Pedersen--Weibel as in Definition \refD{PWcats} with the associated $K$-theory spectrum $K (X \times Y, R)$.
We construct another category associated to the product $X \times Y$.

Starting with this section, we will assume that each metric space we consider has a single metric component.  In other words, the metric function has only finite values.

\begin{DefRef}{FBCcat}
The new category has the same objects as $\mathcal{C} (X \times Y, R)$ but a weaker control condition on the morphisms.
For a choice of a base point $x_0$ in $X$, any function $f \colon [0, +\infty) \to [0, +\infty)$, and a real number $D \ge 0$, define
\[
N (D,f) (x,y) = x[D] \times y[f(d(x,x_0))],
\]
the $(D,f)$-neighborhood of $(x,y)$ in $X \times Y$.

A homomorphism $\phi \colon F \to G$ is called \textit{fibrewise} $(D,f)$-\textit{bounded} for some choice of $D \ge 0$ and a non-decreasing function $f \colon [0, +\infty) \to [0, +\infty)$, or simply \textit{fibrewise bounded}, if the components $F_{(x,y)} \to G_{(x',y')}$ are zero maps for $(x',y')$ outside of the $(D,f)$-neighborhood of $(x,y)$.  
\end{DefRef}

\begin{LemRef}{OKLIAS}
	The notion of a fibrewise bounded homomorphism is independent of the choice of $x_0$ in $X$.
\end{LemRef}

\begin{proof} 
	It suffices to show that if $\phi$ is fibrewise $(D,f)$-bounded with respect to a choice of a base point $x_0$, and if given a different choice $x'_0$ of a base point, then $\phi$ is fibrewise $(D',f')$-bounded for some choice of the parameters $(D',f')$. 
	We choose $D' = D$ and define $f' (t) = f (t + d(x_0, x'_0))$.  From the triangle inequality, $d(x,x_0) \le d(x,x'_0) + d(x_0, x'_0)$, so \[ f' (d(x,x'_0)) = f (d(x,x'_0) + d(x_0, x'_0)) \ge f(d(x,x_0)).\]  This shows that $x[D'] \times y[f'(d(x,x'_0))]$ always contains $x[D] \times y[f(d(x,x_0))]$.  
\end{proof}

We define the \textit{category of geometric modules over the product $X \times Y$ with fibrewise control over $X$} as the category of usual geometric $R$-modules over the product metric space and fibrewise bounded homomorphisms.

\begin{NotRef}{KINYBTR}
	The notation for this category which emphasizes the special role of the factor $X$ is $\mathcal{C}_X (Y,R)$.  This is very much in line with the original notation $\mathcal{C}_X (R)$ used by Pedersen and Weibel for the bounded category of $R$-modules over $X$ and specifically $\mathcal{C}_i (R)$ for the theory over $\mathbb{R}^i$, see for example 
	\cite[Remark 1.2.3]{ePcW:85}.
	To streamline the notation even further, we often omit the the ring $R$ from the notation $\mathcal{C}_X (Y)$ as it usually plays the role of a dummy variable.
\end{NotRef}

\begin{DefRef}{FtF}
The connective \textit{bounded $K$-theory of geometric $R$-modules over the product $X \times Y$ with fibrewise control over $X$} is the spectrum $K_X (Y)$ associated to the additive category $\mathcal{C}_X (Y)$. 
\end{DefRef}

It follows from Example 1.2.2 of Pedersen-Weibel \cite{ePcW:85} that in general the fibred bounded category $\mathcal{C}_X (Y)$ is not
isomorphic to $\mathcal{C} (X \times Y, R)$.
The proper generality of that work, as explained in \cite{ePcW:85,ePcW:89}, starts with a general additive category $\mathcal{A}$ embedded in a cocomplete additive category, generalizing the setting of free finitely generated $R$-modules as a subcategory of all free $R$-modules.
All of the results of Pedersen-Weibel hold for $\mathcal{C} (X, \mathcal{A})$.

In these terms, the category $\mathcal{C}_X (Y)$ can be seen isomorphic to the category $\mathcal{C} (X, \mathcal{A})$, where $\mathcal{A} = \mathcal{C} (Y, R)$.  Unfortunately this concise description is not sufficient for defining features such as the necessary choice of Karoubi filtrations and other details in forthcoming arguments.  

The difference between $\mathcal{C}_X (Y)$ and $\mathcal{C} (X \times Y, R)$ is made to disappear in \cite{ePcW:85}
by making $\mathcal{C} (Y,R)$ ``remember the filtration'' of morphisms when viewed
as a filtered additive category with
$\Hom_D (F,G)$ consisting of all morphisms $\phi \in \Hom (F,G)$ which are bounded by $D$.
Identifying a small category with its set of morphisms, one can
think of the bounded category as
\[
  \mathcal{C} (Y,R) = \colim{D \in \mathbb{R}} \mathcal{C}_D (Y,R),
\]
where $\mathcal{C}_D (Y,R) = \Hom_D (\mathcal{C} (Y,R))$ is the
collection of all $\Hom_D (F,G)$. 
This interpretation gives an exact embedding 
\[
\mathcal{C} (X \times Y, R)
= \colim{D \in \mathbb{R}} \mathcal{C} (X,
\mathcal{C}_D (Y,R)) \xrightarrow{\ \iota \ }  \mathcal{C} (X, \colim{D \in \mathbb{R}} \mathcal{C}_D(Y, R)) =  \mathcal{C} (X, \mathcal{C}(Y, R)),
\]
which induces a map of $K$-theory spectra $K(\iota) \colon K (X \times Y, R) \to K_X (Y)$.

We want to develop some results for a variety of categories with fibred bounded control where the Karoubi filtration techniques suffice.

\begin{NotRef}{JUODQ}
Let
\begin{align}
\mathcal{C}_k = \ &\mathcal{C}_{X} (Y \times \mathbb{R}^k), \notag \\
\mathcal{C}^{+}_k = \ &\mathcal{C}_{X} (Y \times \mathbb{R}^{k-1} \times [0,+\infty)), \notag \\
\mathcal{C}^{-}_k = \ &\mathcal{C}_{X} (Y \times \mathbb{R}^{k-1} \times (-\infty,0]). \notag
\end{align}
We will also use the notation
\begin{align}
\mathcal{C}^{<+}_k =& \colim{D \ge 0} \mathcal{C}_{X} (Y \times \mathbb{R}^{k-1} \times [-D,+\infty)), \notag \\
\mathcal{C}^{<-}_k =& \colim{D \ge 0} \mathcal{C}_{X} (Y \times \mathbb{R}^{k-1} \times (-\infty,D]), \notag \\
\mathcal{C}^{<0}_k =& \colim{D \ge 0} \mathcal{C}_{X} (Y \times \mathbb{R}^{k-1} \times [-D,D]). \notag
\end{align}
\end{NotRef}
Clearly $\mathcal{C}_k$ is $\mathcal{C}^{<-}_k$-filtered and that $\mathcal{C}^{<+}_k$ is $\mathcal{C}^{<0}_k$-filtered.  There are equivalences of categories $\mathcal{C}^{<0}_k \simeq \mathcal{C}_{k-1}$, $\mathcal{C}^{<-}_k \simeq \mathcal{C}^{-}_k$, and $\mathcal{C}_k/ \mathcal{C}^{<-}_k \simeq \mathcal{C}^{<+}_k/ \mathcal{C}^{<0}_k$, as explained at the end of section \refS{BETRGCV}.  By Theorem \refT{CPF}, the commutative diagram
\[
\begin{CD}
K ((\mathcal{C}^{<0}_k)^{\wedge K}) @>>> K (\mathcal{C}^{<+}_k) @>>> K (\mathcal{C}^{<+}_k/ \mathcal{C}^{<0}_k) \\
@VVV @VVV @VV{\cong}V \\
K ((\mathcal{C}^{<-}_k)^{\wedge K}) @>>> K (\mathcal{C}_k) @>>> K (\mathcal{C}_k/ \mathcal{C}^{<-}_k)
\end{CD}
\]
where all maps are induced by inclusions on objects, is in fact a map of homotopy fibrations.  The categories $\mathcal{C}^{<+}_k$ and $\mathcal{C}^{<-}_k$ are flasque, that is, possess an endofunctor $\Sh$ such that $\Sh (F) \cong F \oplus \Sh (F)$, which can be seen by the usual Eilenberg swindle argument.
Therefore $K (\mathcal{C}^{<+}_k)$ and $K (\mathcal{C}^{<-}_k)$ are contractible by the
Additivity Theorem, cf.\ Pedersen--Weibel \cite{ePcW:85}.   This gives a map $K (\mathcal{C}_{k-1}) \to \Omega  K (\mathcal{C}_k)$ which induces isomorphisms of $K$-groups in positive dimensions.

\begin{DefRef}{POIYSNHU}
The \textit{nonconnective fibred bounded $K$-theory} is the spectrum
\[
\KncX (Y) \overset{ \text{def} }{=} \hocolim{k>0}
\Omega^{k} K (\mathcal{C}_k).
\]
\end{DefRef}

If $Y$ is the single point space then the delooping $\KncX (\point)$ is clearly equivalent to
the nonconnective delooping $\Knc (X,R)$ of Pedersen--Weibel via the map
$K(\iota) \colon \Knc (X \times \point, R) \rightarrow \KncX (\point)$.

\begin{RemRef}{Other}
From general principles, the delooping here is formally equivalent to the one that is obtained by using a different model $\mathcal{C}_{X \times \mathbb{R}^k} (Y)$ instead of $\mathcal{C}_{X} (Y \times \mathbb{R}^k)$.  In the latter approach one simply uses the usual bounded excision for the decomposition of $X \times \mathbb{R}^k$ as the union of its subsets 
$X \times \mathbb{R}^{k-1} \times [0,+\infty)$ and 
$X \times \mathbb{R}^{k-1} \times (-\infty,0]$.  The details of our delooping in Definition \refD{POIYSNHU} are important because they are compatible with the upcoming required fibrewise deloopings.
\end{RemRef}

Bounded excision theorems of section \refS{BETRGCV} are easily adapted to the fibred setting of 
$\KncX (Y)$.
We also include concise statements made in terms of coverings by ``coarse families'' in the language introduced in \cite[section 4.3]{gCbG:19}. 

Suppose $Y_1$
and $Y_2$ are mutually antithetic subsets of a proper metric space $Y$, and $Y = Y_1
\cup Y_2$.

\begin{NotRef}{ZXCVRE}
	Let us use the general notation $\mathcal{C}_X
(Y)_{<C}$, for a subset $C$ of $Y$, to denote the full subcategory of $\mathcal{C}_X (Y)$ on the objects with supports $Z \subset X \times Y$ which are coarsely equivalent to the subset $X \times C$.
	
	We will also use the shorthand notation $\mathcal{C} = \mathcal{C}_X (Y)$, $\mathcal{C}_i= \mathcal{C}_X
(Y)_{<Y_i}$ for indices $i=1$ or $2$, and $\mathcal{C}_{12}$ for the intersection
$\mathcal{C}_1 \cap \mathcal{C}_2$.
\end{NotRef}

There is a commutative diagram
\[
\begin{CD}
K (\mathcal{C}_{12}) @>>> K (\mathcal{C}_1) @>>> K ({\mathcal{C}_1}/{\mathcal{C}_{12}}) \\
@VVV @VVV @VV{K(I)}V \\
K (\mathcal{C}_2) @>>> K (\mathcal{C}) @>>> K ({\mathcal{C}}/{\mathcal{C}_2})
\end{CD} \tag{$\natural$}
\]
where the rows are homotopy fibrations from Theorem \refT{CPF} and $I \colon {\mathcal{C}_1}/{\mathcal{C}_{12}} \to {\mathcal{C}}/{\mathcal{C}_2}$ is the functor induced from the exact inclusion $I \colon \mathcal{C}_1 \to \mathcal{C}$.
It clear that $I$ is again an additive equivalence of categories.

The subcategory $\mathcal{C}_{X} (Y \times \mathbb{R}^k)_{<C \times \mathbb{R}^k}$ is evidently a Karoubi subcategory of $\mathcal{C}_{X} (Y \times \mathbb{R}^k)$ for any choice of the subset $C \subset Y$.
We define
\[
\KncX (Y)_{<C} \overset{ \text{def} }{=} \hocolim{k>0}
\Omega^{k} {K}_{X} (Y \times \mathbb{R}^k)_{<C \times \mathbb{R}^k}.
\]
Using the methods above, one easily obtains the weak equivalence
\begin{gather*}
\KncX (Y)_{<C} \simeq \KncX (C).
\end{gather*}
We also define
\[
\KncX (Y)_{<Y_1, Y_2} \overset{ \text{def} }{=} \hocolim{k>0}
\Omega^{k} {K}_{X} (Y \times \mathbb{R}^k)_{<Y_1 \times \mathbb{R}^k, \, Y_2 \times \mathbb{R}^k}.
\]

We are now able to closely follow the sketch of proof of Theorem \refT{BddExc} to get the following result. 

\begin{ThmRefName}{Exc2}{Fibrewise Bounded Excision}
Suppose $Y_1$ and $Y_2$ are subsets of a metric space $Y$, and $Y = Y_1 \cup Y_2$.
There is a homotopy pushout diagram of spectra
\[
\begin{CD}
\KncX (Y)_{<Y_1,Y_2} @>>> \KncX (Y)_{<Y_1} \\
@VVV @VVV \\
\KncX (Y)_{<Y_2} @>>> \KncX (Y)
\end{CD}
\]
where the maps are induced from the exact inclusions.
If $Y_1$ and $Y_2$ are mutually antithetic subsets of $Y$,
there is a homotopy pushout
\[
\begin{CD}
\KncX (Y_1 \cap Y_2) @>>> \KncX (Y_1) \\
@VVV @VVV \\
\KncX (Y_2) @>>> \KncX (Y)
\end{CD}
\]
\end{ThmRefName}

Using the terminology from \cite[section 4.3]{gCbG:19}, suppose $\mathcal{U}$ is a finite coarse covering of $Y$ of cardinality $s$, closed under coarse intersections, and such that the family of all subsets $U$ in $\mathcal{U}$ is pairwise coarsely antithetic.
The coarsely saturated families which are members of $\mathcal{U}$ are partially ordered by inclusion.
In fact, the union of the families forms the set $\mathcal{A}$ closed under intersections.

\begin{CorRef}{gnaloi}
We define the homotopy pushout
\[
\KncX (Y; \mathcal{U}) = \hocolim{A \in \mathcal{A} \in \mathcal{U}} \KncX (Y)_{<A}.
\]
Then the map \[ \KncX (Y; \mathcal{U})  \longrightarrow 
\KncX (Y) \]
induced by inclusions $\mathcal{C}_{X} (Y)_{<U} \to \mathcal{C}_{X} (Y)$
is a weak equivalence.
\end{CorRef}

\begin{proof}
Apply Theorem \refT{Exc2} inductively to the sets in $\mathcal{U}$.
\end{proof}

There is also a relative version of fibred $K$-theory and the corresponding Fibrewise Excision Theorem.

\begin{DefRef}{GEPCWpr}
Let $Y' \in \mathcal{A}$ as part of a coarse covering $\mathcal{U}$ of $Y$.  Let $\mathcal{C} = \mathcal{C}_X (Y)$ and $\mathcal{Y}' = \mathcal{C}_X (Y)_{<Y'}$.  The category $\mathcal{C}_X (Y,Y')$ is the quotient category $\mathcal{C}/\mathcal{Y}'$.

It is now straightforward to define
\[
\KncX (Y,Y') {=} \hocolim{k>0}
\Omega^{k} {K}_{X} (Y \times \mathbb{R}^k, Y' \times \mathbb{R}^k),
\]
\[
\KncX (Y,Y')_{<C} {=} \hocolim{k>0}
\Omega^{k} {K}_{X} (Y \times \mathbb{R}^k, Y' \times \mathbb{R}^k)_{<C \times \mathbb{R}^k},
\]
and
\[
\KncX (Y,Y')_{<C_1, C_2} {=} \hocolim{k>0}
\Omega^{k} {K}_{X} (Y \times \mathbb{R}^k, Y' \times \mathbb{R}^k)_{<C_1 \times \mathbb{R}^k, \, C_2 \times \mathbb{R}^k}.
\]
\end{DefRef}

The theory developed in this section is easily relativized to give the following excision theorems.

\begin{ThmRefName}{ExRel}{Relative Fibrewise Excision}
If $Y$ is the union of two subsets $Y_1$ and $Y_2$,
there is a homotopy pushout of spectra
\[
\begin{CD}
\KncX (Y,Y')_{<Y_1,Y_2} @>>> \KncX (Y,Y')_{<Y_1} \\
@VVV @VVV \\
\KncX (Y,Y')_{<Y_2} @>>> \KncX (Y,Y')
\end{CD}
\]
where the maps are induced from the exact inclusions.
In fact, if $Y$ is the union of two mutually antithetic subsets $Y_1$ and $Y_2$, and $Y'$ is antithetic to both $Y_1$ and $Y_2$,
there is a homotopy pushout
\[
\begin{CD}
\KncX (Y_1 \cap Y_2, Y_1 \cap Y_2 \cap Y') @>>> \KncX (Y_1, Y_1 \cap Y') \\
@VVV @VVV \\
\KncX (Y_2, Y_2 \cap Y') @>>> \KncX (Y,Y')
\end{CD}
\]
More generally, suppose $\mathcal{U}$ is a finite coarse covering of $Y$ by mutually antithetic subsets, all antithetic to $Y'$.
We define the homotopy pushout
\[
\KncX (Y,Y'; \mathcal{U}) = \hocolim{A \in \mathcal{A} \in \mathcal{U}} \KncX (Y,Y')_{<A}.
\]
Then there is a weak equivalence
\[
\KncX (Y,Y'; \mathcal{U}) \simeq \KncX (Y,Y').
\]
\end{ThmRefName}

\begin{ThmRef}{ExcExcPrep}
Given a subset $U$ of $Y'$, there is an equivalence
\[
\KncX (Y,Y') \simeq \KncX (Y-U,Y'-U).
\]
\end{ThmRef}

\begin{proof}
Consider the setting of Theorem \refT{Exc2} with $Y_1 = Y-U$ and $Y_2=Y'$. The map
\[
{\mathcal{C}_X (Y)_{<(Y-U)}}/{\mathcal{C}_X (Y)_{<(Y-U)} \cap \mathcal{C}_X (Y)_{<Y'}} \longrightarrow {\mathcal{C}_X (Y)}/{\mathcal{C}_X (Y)_{<Y'}}
\]
is an equivalence of categories as in the sketch of proof of Theorem \refT{BddExc}, and so induces a weak equivalence of $K$-theory spectra.
Since $U$ is a subset of $Y'$,
\[
{\mathcal{C}_X (Y)_{<(Y-U)} \cap \mathcal{C}_X (Y)_{<Y'}} = \mathcal{C}_X (Y)_{<(Y'-U)}.
\]
Now the maps of quotients
\[
{\mathcal{C}_X (Y)}/{\mathcal{C}_X (Y')}  \longrightarrow  {\mathcal{C}_X (Y)}/{\mathcal{C}_X (Y)_{<Y'}},
\]
\[
  {\mathcal{C}_X (Y-U)}/{\mathcal{C}_X (Y'-U)}
  \longrightarrow  
  {\mathcal{C}_X (Y)_{<(Y-U)}}/{\mathcal{C}_X (Y)_{<(Y'-U)}}
\]
also induce weak equivalences.
Their composition gives the required equivalence.
\end{proof}

\SecRef{Equivariant parametrized \textit{K}-theory}{JKASEQT}

The classical situation is a proper metric space $M$ with a free, properly discontinuous left $\Gamma$-action by isometries.  
In this case,
there is clearly a natural action of $\Gamma$ on the geometric modules $\mathcal{C} (M)$ and
therefore on $\Knc (M)$.  Formally, this action is also free.
A different equivariant bounded $K$-theory with useful fixed points is constructed as follows.

\begin{DefRef}{EqCatDef}
Let $\EGamma$ be the category with the object set $\Gamma$ and the
unique morphism $\mu \colon \gamma_1 \to \gamma_2$ for any pair
$\gamma_1$, $\gamma_2 \in \Gamma$. There is a left action by $\Gamma$
on $\EGamma$ induced from the left multiplication in $\Gamma$.
\end{DefRef}

If $\mathcal{C}$ is a small category with a left $\Gamma$-action,
then the functor category
$\Fun(\EGamma,\mathcal{C})$ is a category with the
left $\Gamma$-action given on functors $\theta \colon \EGamma \to \mathcal{C}$ by
$\gamma(\theta)(\gamma')=\gamma \theta (\gamma^{-1} \gamma')$
and
$\gamma(\theta)(\mu)=\gamma \theta (\gamma^{-1} \mu)$.
It is
nonequivariantly equivalent to $\mathcal{C}$.

\begin{DefRef}{LaxLimit}
The subcategory of equivariant functors and
equivariant natural transformations in $\Fun(\EGamma,\mathcal{C})$
is the fixed subcategory $\Fun(\EGamma,\mathcal{C})^{\Gamma}$ known as
the \textit{lax limit of the action}, introduced and called so by Thomason in \cite{rT:83}.

Thomason also gave the following explicit description of the lax limit.
An object $\theta$ of
$\Fun(\EGamma,\mathcal{C})^{\Gamma}$ is determined by its value $F = \theta (e) \in \mathcal{C}$ on the identity element $e$ of $\Gamma$ together with a function $\psi$ defined on $\Gamma$ with $\psi (\gamma) \in  \Hom (F, \gamma F )$, subject to the condition $\psi(e) = \id$ and the cocycle identity $\psi (\gamma_1 \gamma_2) = \gamma_1 \psi(\gamma_2) \circ  \psi (\gamma_1)$.  These conditions imply that $\psi (\gamma)$ is always an
isomorphism.  

Identifying the objects of $\Fun(\EGamma,\mathcal{C})^{\Gamma}$ with pairs
$(F,\psi)$ as above, a natural transformation $\theta \to \theta'$ of functors in $\Fun(\EGamma,\mathcal{C})^{\Gamma}$ is completely determined by a morphism $\phi \colon F \to  F'$ in
$\mathcal{C}$ such that the squares 
\[
\begin{CD}
F @>{\psi (\gamma)}>> \gamma F \\
@V{\phi}VV @VV{\gamma \phi}V \\
F' @>{\psi' (\gamma)}>> \gamma F'
\end{CD}
\]
commute for all $\gamma \in \Gamma$.
\end{DefRef}

Specializing this definition to the case of $\mathcal{C} = \mathcal{C} (M)$, we will use the notation
$\mathcal{C}^{\Gamma} (M)$ for the category $\Fun(\EGamma,\mathcal{C} (M))$ with the $\Gamma$-action defined above.
It turns out that the following adjustment makes the fixed points better behaved.
Notice that
$\mathcal{C} (M)$ contains the family of isomorphisms $\phi$ such that $\phi$ and $\phi^{-1}$ are bounded by $0$.
We will express this property by saying that the filtration of $\phi$ is $0$ and writing $\fil (\phi) =0$.
The full
subcategory of functors $\theta \colon \EGamma \to \mathcal{C} (M)$ such that $\fil \theta
(\mu) = 0$, for all $\mu$, is invariant under the $\Gamma$-action.

\begin{DefRef}{PrKeq}
Let $\mathcal{C}^{\Gamma,0} (M)$
be the full subcategory of $\mathcal{C}^{\Gamma} (M)$ on the functors sending all
morphisms in $\EGamma$ to filtration $0$ homomorphisms.
Then $\mathcal{C}^{\Gamma,0} (M)^{\Gamma}$  is the full subcategory of $\mathcal{C}^{\Gamma} (M)^{\Gamma}$ on the objects
$(F, \psi)$ with $\fil \psi (\gamma) = 0$ for all $\gamma \in \Gamma$.
\end{DefRef}

We define
$K^{\Gamma,0} (M)$ to be the nonconnective delooping of the $K$-theory of the symmetric monoidal category $\mathcal{C}^{\Gamma,0} (M)$.

It is shown in section VI of \cite{gC:95} that the fixed points of the induced $\Gamma$-action on $K^{\Gamma,0} (M)$ is the nonconnective delooping of the $K$-theory of
$\mathcal{C}^{\Gamma,0} (M)^{\Gamma}$.
One also knows the following.

\begin{ThmRefName}{VVV}{Corollary VI.8 of \cite{gC:95}}
If $M$ is a proper metric space and $\Gamma$ acts on $M$ freely, properly discontinuously, cocompactly by isometries, there are weak equivalences
$K^{\Gamma,0} (M)^{\Gamma} \simeq 
\Knc (M/\Gamma, R[\Gamma]) \simeq \Knc (R[\Gamma])$.
\end{ThmRefName}

This theorem applies in two specific cases of interest to us.
One case is of the
fundamental group $\Gamma$ of a closed aspherical manifold acting on the universal cover by
covering transformations.
The second is of $\Gamma$ acting on itself, as a word-length metric space, by left
multiplication.

\begin{RemRef}{NotSame}
The theory $K^{\Gamma,0} (M)^{\Gamma}$ may very well
differ from $K^{\Gamma} (M)^{\Gamma}$.
According to Theorem \refT{VVV}, the fixed point category
of the correct bounded equivariant theory is, for example, the category of free
$R[\Gamma]$-modules when $M = \Gamma$ with the word metric relative to a chosen generating set. 
However,
$\mathcal{C}^{\Gamma} (\Gamma,R)^{\Gamma}$ will
include the $R$-module with a single stalk $F_{\gamma} = R$ for all $\gamma$ and equipped with trivial
$\Gamma$-action, which is not free.
\end{RemRef}

We will treat the group $\Gamma$ equipped with a finite generating set $\Omega$ as a metric space with a word-length metric, as explained in Definition \refD{WORD}.  It is well-known that varying $\Omega$ only changes $\Gamma$ to a
coarsely equivalent metric space, cf. \cite[Proposition 1.15]{jR:03} .

\begin{DefRefName}{POI1}{Coarse equivariant theories}
We associate two new equivariant theories on metric spaces with a left $\Gamma$-action.  The theory $K_i^{\Gamma}$ is defined only for metric
spaces $Y$ with actions by isometries, while $K_p^{\Gamma}$ for metric
spaces with coarse actions.
\begin{enumerate}
\item
$k^{\Gamma}_i (Y)$ is defined to be the $K$-theory of
$\mathcal{C}^{\Gamma}_i (Y) = \mathcal{C}^{\Gamma,0}(\Gamma \times Y, R)$,
where $\Gamma$ is regarded as a
word-length metric space with isometric $\Gamma$-action given by left multiplication, and $\Gamma
\times Y$ is given the product metric and the diagonal isometric action.
\item
$k^{\Gamma}_p (Y)$ is defined for any metric space $Y$ equipped with a $\Gamma$-action
by coarse equivalences. It is the $K$-theory spectrum attached to a symmetric
monoidal category $\mathcal{C}^{\Gamma}_p (Y)$ with $\Gamma$-action whose objects are given by functors
\[
\theta \colon \EGamma \longrightarrow \mathcal{C}_{\Gamma} (Y) =
\mathcal{C} (\Gamma, \mathcal{C} (Y,R))
\]
with the additional condition that the morphisms $\theta(\mu)$ are of filtration zero but only as homomorphisms between $R$-modules parametrized over $\Gamma$.
\end{enumerate}
\end{DefRefName}

Now the nonconnective equivariant $K$-theory spectra $K_i^{\Gamma}$ and $K_p^{\Gamma}$ should be the nonconnective deloopings of $k_i^{\Gamma}$ and $k_p^{\Gamma}$.
For example, if we define
\[
\mathcal{C}_{i,k}^{\Gamma} = \mathcal{C}^{\Gamma,0}(\Gamma \times \mathbb{R}^{k} \times Y, R),
\]
where $\Gamma$ acts on the product $\Gamma \times \mathbb{R}^{k} \times Y$ according to $\gamma  (\gamma', x, y) = (\gamma \gamma', x, \gamma(y))$, and
\[
\mathcal{C}_{i,k}^{\Gamma,+} = \mathcal{C}^{\Gamma,0}(\Gamma \times \mathbb{R}^{k-1} \times [0,+\infty) \times Y, R), \ \mathrm{etc.,}
\]
then the delooping construction in Definition \refD{POIYSNHU} can be applied verbatim.
Similarly, one can use the $\Gamma$-action on $\Gamma \times \mathbb{R}^{k}$ given by
$\gamma  (\gamma', x) = (\gamma \gamma', x)$ and define $\mathcal{C}_{p,k}^{\Gamma}$ as the symmetric monoidal category of functors
$\theta \colon \EGamma \rightarrow
\mathcal{C} (\Gamma \times \mathbb{R}^{k}, \mathcal{C} (Y,R))$
such that the morphisms $\theta(f)$ are bounded by $0$ as $R$-linear homomorphisms
over $\Gamma \times \mathbb{R}^{k}$.
There are obvious analogues of the categories $\mathcal{C}_{p,k}^{\Gamma,+}$, etc.
If $\ast$ is either of the subscripts $i$ or $p$, and $Y$ is equipped with actions by $\Gamma$ via respectively isometries or coarse equivalences, we obtain equivariant maps
\[
K (\mathcal{C}^{\Gamma}_{\ast, k-1} (Y)) \longrightarrow \Omega K (\mathcal{C}^{\Gamma}_{\ast,k} (Y)).
\]

\begin{DefRef}{EQIOMN}
Let $\ast$ be either of the subscripts $i$ or $p$.
We define
\[
k^{\,\Gamma}_{\ast,k} (Y) = K (\mathcal{C}^{\Gamma}_{\ast,k} (Y))
\]
and the nonconnective equivariant spectra
\[
K^{\Gamma}_{\ast} (Y) \overset{ \text{def} }{=} \hocolim{k>0}
\Omega^{k} k^{\,\Gamma}_{\ast,k} (Y).
\]
The same construction gives for the fixed points
\[
K^{\Gamma}_{\ast} (Y)^{\Gamma} = \hocolim{k>0}
\Omega^{k} k^{\,\Gamma}_{\ast,k} (Y)^{\Gamma}.
\]

Given left $\Gamma$-actions by isometries on metric spaces $X$ and $Y$, there are evident diagonal actions induced on the categories $\mathcal{C} (X \times Y, R)$ and $\mathcal{C}
(X, \mathcal{C} (Y,R))$. The equivariant embedding induces
the equivariant functor
\[
i^{\Gamma} \colon \mathcal{C}^{\Gamma} (X \times Y, R) \longrightarrow
\mathcal{C}^{\Gamma}
(X, \mathcal{C} (Y,R)).
\]
Taking $X = \Gamma$, there results a natural transformation $K^{\Gamma}_i (Y) \to  K^{\Gamma}_p (Y)$.
\end{DefRef}

One basic relation between the two equivariant fibrewise theories is through the observation that
in both cases, when $Y$ is a single point space, $\mathcal{C}^{\Gamma}_i (\point)$  and $\mathcal{C}^{\Gamma}_p (\point)$ can be identified
with $\mathcal{C}^{\Gamma,0} (\Gamma, R)$.

Another fact is a special instructive property of the theory $K_p^{\Gamma}$.

\begin{ThmRef}{GpDef}
Suppose $\Gamma$ acts on a metric space $Y$ by bounded coarse equivalences.  Let $Y_0$ be the same metric space but with $\Gamma$ acting trivially by the identity.  Then there is a weak equivalence
\[
\zeta \colon
K_p^{\Gamma} (Y)^{\Gamma} \xrightarrow{ \ \simeq \ } K_p^{\Gamma} (Y_0)^{\Gamma}.
\]
\end{ThmRef}

\begin{proof}
The category $\mathcal{C}_p^{\Gamma} (Y)$ has the left action by $\Gamma$ induced from the diagonal action on $\Gamma \times Y$.
Recall that an object of $\mathcal{C}_p^{\Gamma} (Y)^{\Gamma}$ is determined by an object
$F$ of $\mathcal{C}_{\Gamma} (Y)$ and isomorphisms $\psi (\gamma) \colon F \to {\gamma} F$ which are of filtration $0$ when projected to $\Gamma$.
Given two objects $(F, \{\psi (\gamma) \})$ and $(G, \{\phi (\gamma) \})$, a morphism $\lambda \colon (F, \{\psi (\gamma) \}) \to (G, \{\phi (\gamma) \})$ is given by a morphism $\lambda \colon F \to G$ in $\mathcal{C}_{\Gamma} (Y)$ such that the collection of morphisms $\gamma \lambda \colon {\gamma} F \to {\gamma} G$ satisfies
\[
\phi (\gamma) \circ \lambda =  \gamma \lambda \circ \psi (\gamma)
\]
for all $\gamma$ in $\Gamma$.
Given $(F, \{\psi (\gamma) \})$, define $(F_0, \{\psi_0 (\gamma) \})$ by $F_0 = F$ and $\psi_0 (\gamma)= \id_{F}$ for all $\gamma \in \Gamma$.  Then $\psi (\gamma)^{-1}$ give a natural isomorphism $Z_F$ from $F$ to $F_0$ and induce
an equivalence
\[
\zeta \colon K_p^{\Gamma} (Y)^{\Gamma} \, \simeq \,  K_p^{\Gamma} (Y_0)^{\Gamma}.
\]
Of course, the bound for the isomorphism $\psi (\gamma)^{-1}$ can vary with $\gamma$.
\end{proof}

Notice that this equivalence exists only in the theory $K_p^{\Gamma}$ and not in $K_i^{\Gamma}$ because the homomorphism $\lambda$ we used in the proof is not a bounded homomorphism and is a morphism only in $\mathcal{C}_p^{\Gamma}$.

Throughout this section we will fix one left action of $\Gamma$ on $Y$ by bounded coarse equivalences in the sense of Definition \refD{leftbdd}.  This means that for each element $\gamma$ the corresponding self-equivalence of $Y$ is a bounded map, but the bound is allowed to vary with $\gamma$.  

We will state results only for the equivariant theory $K_p^{\Gamma}$ where they are the most useful. Some but not all of the statements are also true in $K_i^{\Gamma}$.  However, it is easy to see that all of the statements we make are true for $K_i^{\Gamma}$ if the action is the trivial action by the identity for all $\gamma$.

\begin{DefRef}{KQBdd}
Using the conventions from Notation \refN{ZXCVRE}, let
$\mathcal{C}^{\Gamma}_{p} (Y)_{<Y'}$ be the full subcategory of $\mathcal{C}^{\Gamma}_{p} (Y)$ on objects $\theta$ such that the support of each $\theta (\gamma)$ is contained in a subset coarsely equivalent to $\Gamma \times Y'$.
From viewing this as a parametrized version of $\mathcal{C} (M)_{< U}$-filtrations of $\mathcal{C}(M)$ in section \refS{BETRGCV}, it is clear that the subcategories give Karoubi filtrations and therefore Karoubi quotients
$\mathcal{C}^{\Gamma}_{p} (Y,Y')$.
The bounded actions of $\Gamma$ extend to the quotients in each case.
Taking the $K$-theory of the symmetric monoidal categories with $\Gamma$-actions gives $\Gamma$-spectra
$k^{\Gamma}_{p} (Y,Y')$.

One can now construct the parametrized versions of the relative module categories
$\mathcal{C}^{\Gamma}_{p,k} (Y,Y')$,
their $K$-theory spectra $k^{\Gamma}_{p,k} (Y,Y')$, and the resulting deloopings.
Thus we obtain the nonconnective $\Gamma$-spectra
\[
K^{\Gamma}_{p} (Y,Y') \overset{ \text{def} }{=} \hocolim{k>0}
\Omega^{k} k^{\,\Gamma}_{p,k} (Y,Y').
\]
\end{DefRef}

The Fibration Theorem \refT{CPF} can be applied to prove generalizations of bounded excision.  

\begin{DefRef}{THISONE}
The quotient map of categories induces the equivariant map
$K^{\Gamma}_p (Y) \rightarrow
K^{\Gamma}_p (Y,Y')$
and
the map of the fixed points
$K^{\Gamma}_p (Y)^{\Gamma} \rightarrow
K^{\Gamma}_p (Y,Y')^{\Gamma}$.
\end{DefRef}

The assumption that the action of the group on $Y$ is left-bounded is crucial in this definition.  This allows to extend the group action to the filtering subcategory and finally form the Karoubi quotients that give the relative theory. 

\begin{PropRef}{FibPair}
In the action of $\Gamma$ on $Y$ is left-bounded, there is a homotopy fibration
\[
K^{\Gamma}_p ( Y' )^{\Gamma} \longrightarrow
K^{\Gamma}_p ( Y )^{\Gamma} \longrightarrow K^{\Gamma}_p ( Y,Y' )^{\Gamma}.
\]
\end{PropRef}

\begin{proof}
The fibration follows from the fact that $\mathcal{C}^{\Gamma}_p ( Y )^{\Gamma}$ is
$\mathcal{C}^{\Gamma}_p ( Y )_{<Y'}^{\Gamma}$-filtered.
To identify the fibre as $K^{\Gamma}_p ( Y' )^{\Gamma}$, observe that
the inclusion of the subspace $Y'$ induces equivalences of categories
$\mathcal{C} (Y') \simeq \mathcal{C} (Y)_{<Y'}$ and $\mathcal{C}^{\Gamma}_p ( Y' )^{\Gamma} \simeq \mathcal{C}^{\Gamma}_p ( Y )_{<Y'}^{\Gamma}$.
\end{proof}

\begin{ThmRefName}{BDDEXCI25}{Bounded excision}
If $U_1$ and $U_2$ are a coarsely antithetic pair of subsets of $Y$ which form a cover of $Y$, and the action of $\Gamma$ on $Y$ is bounded, then
\[
\xymatrix{
 K^{\Gamma}_p (U_1 \cap U_2)^{\Gamma} \ar[r] \ar[d]
&K^{\Gamma}_p (U_1)^{\Gamma} \ar[d] \\
 K^{\Gamma}_p (U_2)^{\Gamma} \ar[r]
&K^{\Gamma}_p (Y)^{\Gamma}
}
\]
is a homotopy pushout.
\end{ThmRefName}

\begin{proof}
In view of the isomorphism $\mathcal{C}^{\Gamma}_p (U_1, U_1 \cap U_2) \cong \mathcal{C}^{\Gamma}_p (Y, U_2)$,
we have the weak equivalence
$K^{\Gamma}_p (U_1, U_1 \cap U_2)^{\Gamma} \ \simeq \ K^{\Gamma}_p (Y, U_2)^{\Gamma}$.
We have a map of homotopy fibrations
\[
\begin{CD}
K^{\Gamma}_p (U_1 \cap U_2)^{\Gamma} @>>> K^{\Gamma}_p (U_1)^{\Gamma} @>>> K^{\Gamma}_p (U_1, U_1 \cap U_2)^{\Gamma} \\
@VVV @VVV @VV{\simeq}V \\
K^{\Gamma}_p (U_2)^{\Gamma} @>>> K^{\Gamma}_p (Y)^{\Gamma} @>>> K^{\Gamma}_p (Y, U_2)^{\Gamma}
\end{CD}
\]
which gives the homotopy pushout.
\end{proof}

There are associated relative versions of the excision theorems.

\begin{DefRef}{CoarseInvFun}
Generalizing Definition \refD{THISONE}, if $Y''$ is another coarsely invariant subset of $Y$, then the intersection $Y'' \cap Y'$ is coarsely invariant in both $Y$ and $Y'$, there is an equivariant map
$K^{\Gamma}_p (Y'', Y'' \cap Y')
\rightarrow
K^{\Gamma}_p (Y, Y')$
and on the fixed points
$K^{\Gamma}_p (Y'', Y'' \cap Y')^{\Gamma}
\rightarrow
K^{\Gamma}_p (Y, Y')^{\Gamma}$.
\end{DefRef}

\begin{ThmRefName}{BDDEXCI2}{Relative bounded excision}
Suppose $U_1$, $U_2$, and $Y'$ are three pairwise coarsely antithetic subsets of $Y$
such that
$U_1$ and $U_2$ form a cover of $Y$.
Assuming a bounded action of $\Gamma$ on $Y$,
the commutative square
\[
\xymatrix{
 K^{\Gamma}_{p} (U_1 \cap U_2, Y' \cap U_1 \cap U_2)^{\Gamma} \ar[r] \ar[d]
&K^{\Gamma}_{p} (U_1, Y' \cap U_1)^{\Gamma} \ar[d] \\
 K^{\Gamma}_{p} (U_2, Y' \cap U_2)^{\Gamma} \ar[r]
&K^{\Gamma}_{p} (Y,Y')^{\Gamma} }
\]
induced by inclusions of pairs is
a homotopy pushout.
\end{ThmRefName}

\begin{proof}
This follows from the fact that whenever $C$ is a subset of $Y$ which is coarsely antithetic to $Y'$, the category
$\mathcal{C}^{\Gamma}_p ( Y,Y' )^{\Gamma}$ is
$\mathcal{C}^{\Gamma}_p ( Y,Y' )_{<C}^{\Gamma}$-filtered and $\mathcal{C}^{\Gamma}_p ( Y,Y' )_{<C}^{\Gamma}$ is isomorphic to
$\mathcal{C}^{\Gamma}_p ( C,Y' \cap C )^{\Gamma}$.
The details are left to the reader. The rest of the argument is, up to notational complexity, in the sketch of proof of Theorem \refT{BddExc}.
\end{proof}

One gets the following excision statement as a corollary.

\begin{ThmRef}{ExcExc}
Suppose that a proper subset $U$ of $Y$ is coarsely invariant under the action of $\Gamma$ and that $U$, $Y - U$, and $Y'$ form a coarse covering by pairwise antithetic subsets.
If the action of $\Gamma$ is bounded, then there is a weak equivalence
\[
K_p^{\Gamma} (Y-U,Y'-U)^{\Gamma} \longrightarrow K_p^{\Gamma} (Y,Y')^{\Gamma}.
\]
\end{ThmRef}

It is worth pointing out the connection between the ``core'' of the fibred excision computations and the main interest in the applications, the $K$-theory of the group ring $R[\Gamma]$.

\begin{PropRef}{TPRtriv}
Let $y_0$ be a point in $Y$, then
\[
K^{\Gamma}_{p} (Y)_{<y_0} \simeq K^{\Gamma}_{p} (y_0) = K^{\Gamma,0} (\Gamma)^{\Gamma} \simeq \Knc (R[\Gamma]).
\]
\end{PropRef}

\begin{proof}
	The first equivalence is induced from an isomorphism of bounded categories.  The last equivalence is proved in section VI of \cite{gC:95}.
\end{proof}

\SecRef{Example: absence of Karoubi filtrations}{MExD}

We have seen that the lax limit of an action plays an important role in constructing the fixed point spectrum, including the situation with fibred control.  If one wants to compare the effects of two or more actions in the fibre, the appropriate categorical models for the $K$-theory fixed point spectra are based on a relaxed lax limit construction that we describe in the beginning of this section.   We will then demonstrate that the Karoubi filtration techniques are no longer sufficient to analyze $K$-theory of this additive category.

Let $Y$ be a proper metric space with a left action of a finitely generated group $\Gamma$ by coarse equivalences. 

\begin{DefRef}{WGpre2}
The \textit{fibred lax limit} $\mathcal{C}^{h\Gamma} (Y)$ is a category with objects which are sets of data $( \{ F_{\gamma} \},\{ \psi_{\gamma} \} )$ where
\begin{itemize}
\item $F_{\gamma}$ is an object of $\mathcal{C}_{\Gamma} (Y)$ for each $\gamma$ in $\Gamma$,
\item $\psi_{\gamma}$ is an isomorphism $F_e \to F_{\gamma}$ in $\mathcal{C}_{\Gamma} (Y)$,
\item $\psi_{\gamma}$ has filtration $0$ but only when viewed as a morphism in $\mathcal{C} (\Gamma, \mathcal{C}(Y))$,
\item $\psi_e = \id$,
\item $\psi_{\gamma_1 \gamma_2} = \gamma_1 \psi_{\gamma_2} \circ \psi_{\gamma_1}$ for all $\gamma_1$, $\gamma_2$ in $\Gamma$.
\end{itemize}
As usual, the morphisms $( \{ F_{\gamma} \},\{ \psi_{\gamma} \} ) \to ( \{ F'_{\gamma} \},\{ \psi'_{\gamma} \} )$
are collections of morphisms $ \{ \phi_{\gamma} \} $, where each $\phi_{\gamma}$ is a morphism $F_{\gamma} \to F'_{\gamma}$ in $\mathcal{C}_{\Gamma} (Y)$, such that 
\[
\phi_{\gamma} \circ \psi_{\gamma} = \psi'_{\gamma} \circ \phi_e \tag{$\flat$}
\]
for all $\gamma \in \Gamma$.
Notice that in order to define a morphism $\phi$ it suffices to give one component $\phi_e \colon F_e \to F'_e$; all other components $\phi_{\gamma}$ are determined by virtue of $\psi_{\gamma}$ and $\psi'_{\gamma}$ being isomorphisms in ($\flat$).

The additive structure on $\mathcal{C}^{h\Gamma} (Y)$ is induced from that on $\mathcal{C}_{\Gamma} (Y)$.
This means in particular that the operation $\oplus$ has the property $(F \oplus G)_{\gamma} =F_{\gamma} \oplus G_{\gamma}$. 

For any action $\alpha$ on $Y$ by coarse equivalences, the lax limit
$\mathcal{C}_p^{\Gamma} (Y)^{\Gamma}$ is an additive subcategory of $\mathcal{C}^{h\Gamma} (Y)$.
The embedding $E_{\alpha}$ is realized by sending the object $(F, \psi )$ of $\mathcal{C}_{\Gamma} (Y)^{\Gamma}$ to
$( \{ \alpha_{\gamma} F \}, \{ \psi_{\gamma} \} )$.
On the morphisms, 
$E_{\alpha} (\phi) = \{ {\alpha}_{\gamma} \phi \}$.

The spectrum $k^{h\Gamma} (Y)$ is defined as the Quillen $K$-theory of $\mathcal{C}^{h\Gamma} (Y)$.
The exact embeddings $E_{\alpha}$
induce maps of spectra
$\varepsilon_{\alpha}
\colon
k_p^{\Gamma} (Y)^{\Gamma}
\to 
k^{h\Gamma} (Y)$.
The nonconnective delooping $K^{h\Gamma} (Y)$ of $k^{h\Gamma} (Y)$ is constructed along the lines of the delooping $K^{\Gamma}_{p} (Y)^{\Gamma}$ of $k^{\Gamma}_{p} (Y)^{\Gamma}$ in Definition \refD{EQIOMN}. This allows to induce the maps of nonconnective spectra $K^{\Gamma}_{p} (Y)^{\Gamma} \to K^{h\Gamma} (Y)$.
\end{DefRef} 

We proceed to show an elementary example of $\mathcal{C}^{h\Gamma} (Y)$ that illustrates lack of Karoubi filtrations required for fibred excision theorems for $K^{h\Gamma} (Y)$.

The fibre $Y$ is the integers $\mathbb{Z}$ with the standard metric invariant under translations.
The group in this example is the additive group of the integers, the infinite cyclic group, which is denoted by $C$.  It acts on $\mathbb{Z}$ by translations.  When we want to specify the action by an integer $i$, we will use $t_i$ for the corresponding translation.

The goal of the example is to construct a morphism $\phi \colon F \to G$ in the fibred lax limit $\mathcal{C}^{hC} (\mathbb{Z})$ where the codomain lacks necessary $\oplus$-decompositions that should come as part of a Karoubi filtration.  This will be demonstrated by showing that $f$ fails to have the property in the second bullet point in Definition \refD{Karoubi} for any choice of $\oplus$-decomposition of its codomain.  

\begin{NotRef}{HKASD}
We will write $[a,b]$ to mean the subset $\mathbb{Z} \cap [a,b]$ of $\mathbb{Z}$.
\end{NotRef}	

For the basic fibred excision theorem, one should expect to use the filtering subcategory $\mathcal{A} = \mathcal{C}^{hC} (\mathbb{Z})_{< (-\infty, 0]}$, the full subcategory of $\mathcal{C}^{hC} (\mathbb{Z})$ on objects $F$ which have $F_0$ in $\mathcal{C}^{C}_p (\mathbb{Z})_{< (-\infty, 0]}$.

In order to describe $F$ and $G$, we need to specify the collections of objects $F_n$ and $G_n$ and the structure isomorphisms $\psi^F_n$ and $\psi^G_n$.
Let $G_n = \bigoplus_{C \times \mathbb{Z}} R$ and $F_n = \bigoplus_{i \le j} R$ for all $n$.  Notice that $\alpha_k F_n = F_{n+k}$ for the diagonal action $\alpha_k (i,j) = (i+k,j+k)$ on $C \times \mathbb{Z}$, which is the induced action from the choice of the translation action of $C$ on the fibre $\mathbb{Z}$.
We make $F$ an object of $\mathcal{C}^{C}_p (\mathbb{Z})^C$, which means all $\psi^F_n$ are the identity maps.  The structure of $G$ will be more complicated because it will be properly an object of  $\mathcal{C}^{hC} (\mathbb{Z})$.  The isomorphisms $\psi^G_n \colon G_0 \to G_n$ will be $0$-controlled only in the $C$ factor. We proceed to define  $\psi^G_n$.

\begin{NotRef}{HYYY}
Given an object $O$ of $\mathcal{C}^{C}_p (\mathbb{Z})$, we denote by $O(i,a)$ the direct summand indexed by $(i,a)$ in $C \times \mathbb{Z}$.  Similarly, $O (i, [a,b])$ is the direct sum $O(i,a) \oplus O(i,a+1) \oplus \ldots \oplus O(i,b)$, and $O (i) = O (i, (- \infty,+\infty))$. 

We will use matrices to describe $R$-linear automorphisms $\psi^G_n$ of $G$.  Since each $\psi^G_n$ is $0$-controlled along $C$, it is a direct sum of $R$-linear isomorphisms $\psi^G_{n,i}$ between the components $G_{0} (i) \to  G_{n} (i)$.
This means each $\psi^G_{n,i} \colon \bigoplus_{i \times \mathbb{Z}} R \to \bigoplus_{i \times \mathbb{Z}} R$ can be given by a bi-infinite nonsingular matrix over $R$ in the standard basis, with rows and columns enumerated as usual from left to right and from top down respectively. 
Starting with the infinite identity matrix, we replace the square block $[a,b] \times [a,b]$ with the lower triangular matrix of $1$s on and below the diagonal.  Let us denote this matrix by $L (a,b)$.
\end{NotRef}

 The homomorphism $\psi^G_{n,i}$ is defined by the matrix $L (0,n)$.
 
 This is a nonsingular matrix making $\psi^G_{n}$ an isomorphism bounded by $n+1$ in the fibred lax limit $\mathcal{C}^{hC} (\mathbb{Z})$. 
The effect of the action of $C$ on $\psi^G_{n}$ is given by assigning the matrix $L (i,i+n)$ to $t_i \psi^G_{n}$, written in the same basis as $\psi^G_{n}$.  The cocycle identity is satisfied because of the matrix identity $L (i,i+k) L (0,k) = L(0, i+k)$. 

The important feature of $\psi_{n,i}$ is that the submodule $G_0 (0,0)$ maps onto the direct summand in $G_n (0,[0,n])$ generated by the sum of $n+1$ generators $ 1 \oplus 1 \oplus  \ldots 1$.  This direct summand, however, is not a direct summand of the Pedersen-Weibel object $G_n (0)$ in $\mathcal{C} (\mathbb{Z})$ and so not a direct summand of $G_n$ in $\mathcal{C}^{C}_p (\mathbb{Z})$.  If we now restrict $\psi^G_n$ to $G_0 (0, (-\infty,0])$ then the minimal direct summand of $G_n$ containing the image will be the submodule $G_n (0, (-\infty,n])$.  

We remind the reader that in our case the decompositions $E_{\alpha} \oplus D_{\alpha}$ in the definition of Karoubi filtrations are the direct sum decompositions in $\mathcal{C}^{hC} (\mathbb{Z})$.  On the level of each integer $n$, for all $k$ in $\mathbb{Z}$ this enforces direct sum decompositions $G_n (k) = (E_{\alpha})_n (k) \oplus (D_{\alpha})_n (k)$ with the structure maps in $E_{\alpha}$, in particular, being restrictions of $\psi_n$.

Now the morphism $\phi \colon F \to G$ is given by making $\phi_n$ the restriction of $\psi^G_n$ to the submodule $F_n = \bigoplus_{i \le a} G_n (i,a)$.  On the level of $\phi_0 \colon F_0 \to G_0$, this morphism is specified by the inclusion $F_0 \subset G_0$.   We have the commutative square 
\[
\phi_{n} \circ \psi^F_{n} = \psi^G_{n} \circ \phi_0 \tag{$\flat$}
\]
because $\phi_0$ is the identity, $\psi^F_{n}$ is an inclusion, and $\phi_{n}$ is the restriction of $\psi^G_{n}$ by definition.

We claim that $\phi$ fails to factor through any decomposition of $G$ with $E_{\alpha}$ in $\mathcal{A}$.
Indeed, we saw that for each $n$ the smallest direct summand of $G_n$ containing $\phi_0 G_0 (0) = \psi_n G_0 (0)$ is $G_n (0, (-\infty, n])$.  In fact, since $\psi_n$ is an isomorphism and $(\psi_n)_{\alpha}$ for $E_{\alpha}$ is an isomorphism, $E_{\alpha,0} (0)$ should contain $G_0 (0, (- \infty,n]) = \psi_n^{-1} G_n (0, (- \infty,n])$.  Making $n$ arbitrary in this argument shows that $E_{\alpha}$ should contain the entire submodule $G_0 (0, (- \infty, + \infty))$.  This contradicts the fact that $E_{\alpha}$ is an object of $\mathcal{A} = \mathcal{C}^{hC} (\mathbb{Z})_{< (-\infty, 0]}$.

\begin{RemRef}{mvhfbd}
One may wonder about the options available to repair the rigidity of the $\oplus$-structure in $\mathcal{C}^{hC} (\mathbb{Z})$ in order to accommodate images of maps such as $f$.  This would require one to extend the $\oplus$ operation to direct sums of modules which are not themselves Pedersen-Weibel objects even though their sums might be.  To do this coherently, one has to include the filtered free $R$-modules $F$ which assign to each subset $S$ of ${\mathbb{Z}}$, or more generally a proper metric space $M$, a free $R$-submodule $F(S)$ of $F$.  It may be necessary to require that if $S$ is a subset of $T$ then $F(S)$ is a direct summand of $F(T)$.  But unlike the situation in the Pedersen-Weibel objects, one should not expect a splitting $F(T) = F(S) \oplus F(T \setminus S)$.  This means that one should not expect the $K$-theory of this new additive category to be the bounded $K$-theory of $M$. A viable theory of this kind can be constructed and, in fact, does have a useful analogue in $L$-theory.  
                                
Another resolution with good excision properties can be built as the fibred bounded $G$-theory.  This theory was developed in \cite{gCbG:19} and
and possesses not only localization and excision properties of bounded $K$-theory but also localization and excision properties that fail in $K$-theory, as long as the coefficient ring $R$ is Noetherian.  In substantial ways this development mimics the classical relationship between $K$-theory and $G$-theory of Noetherian rings. 
\end{RemRef}

\SecRef{Bounded actions and a cone construction}{PREPSEC}

\begin{DefRef}{OrbitMet}
Let $Z$ be any metric space with a free left $\Gamma$-action by isometries.
We assume that the action is properly discontinuous, that is, for any pair of points $z$ and $z'$,
the infimum over $\gamma \in \Gamma$ of the distances $d(z, \gamma z')$ is attained.
Then we define the orbit space metric on $\Gamma \backslash Z$ by
\[
d_{\Gamma \backslash Z} ([z], [z']) = \inf_{\gamma \in \Gamma} d(z, \gamma z').
\]
\end{DefRef}

\begin{LemRef}{JUQASW}
$d_{\Gamma \backslash Z}$ is a metric on $\Gamma \backslash Z$.
\end{LemRef}

\begin{proof}
It is well-known that $d_{\Gamma \backslash Z}$ is a pseudometric.  The fact that $\Gamma$ acts by isometries makes it a metric.
The triangle inequality follows directly from the triangle inequality for $d$.
Symmetry follows from $d (z, \gamma z') = d(\gamma^{-1} z, z') = d(z', \gamma^{-1} z)$.
Finally, $d_{\Gamma \backslash Z} ([z], [z']) = 0$ gives $d (z, \gamma z') = 0$ for some $\gamma \in \Gamma$, so $d (\gamma' z, \gamma' \gamma z') = 0$ for all $\gamma' \in \Gamma$, and so $[z] = [z']$.
\end{proof}

Now suppose $X$ is some metric space with a left $\Gamma$-action by isometries.

\begin{DefRef}{Xbdd}
Define
\[
X^{bdd} = X \times_{\Gamma} \Gamma
\]
where
the right-hand copy of $\Gamma$ denotes $\Gamma$ regarded as a metric space with the word-length metric associated to a finite generating set,
the group $\Gamma$ acts by isometries on the metric space $\Gamma$ via left multiplication,
and $X \times_{\Gamma} \Gamma$ denotes the orbit metric space associated to the diagonal left $\Gamma$-action
on $X \times \Gamma$.
We will denote the orbit metric by $d^{bdd}$.
\end{DefRef}

A natural left action of $\Gamma$ on $X^{bdd}$ is given by the formula
$\gamma [x, \gamma'] = [x, \gamma' \gamma^{-1}]$.  Then $\gamma [x,e] = [x, \gamma^{-1}] = [\gamma x, e]$.

Suppose $\Gamma$ is a finitely generated group with a chosen finite generating set which determines a word metric $d_{\Gamma}$ and the corresponding norm $\vert \gamma \vert = d_{\Gamma} (e, \gamma)$.

\begin{LemRef}{JUQASW2}
If the left action of  $\Gamma$ on a metric space $X$ is bounded in the sense of Definition \refD{leftbdd}, and $B \colon \Gamma \to [0, \infty)$ is a function such that $d(x, \gamma x) \le B ({\gamma})$ for all $\gamma \in \Gamma$ and $x \in X$, then there is a real function $B_{\ast} \colon [0, \infty) \to [0, \infty)$ such that $\vert \gamma \vert \le s$ implies $B({\gamma}) \le B_{\ast} (s)$. 
\end{LemRef}

\begin{proof}
One simply takes $B_{\ast} (s) = \max \{ B({\gamma}) \ \mathrm{such} \ \mathrm{that} \  \vert \gamma \vert \le s \}$.
\end{proof}

\begin{PropRef}{ECEjkiu}
The natural action of $\Gamma$ on $X^{bdd}$ is bounded.
\end{PropRef}

\begin{proof}
We choose $B_{\gamma} = \vert \gamma \vert$.
Now
\begin{align}
d^{bdd} ([x,e], [\gamma x,e]) &= \inf_{\gamma' \in \Gamma} d^{\times} ((x,e), \gamma' (\gamma x,e)) \notag \\
&\le d^{\times} ((x,e), \gamma^{-1} (\gamma x,e)) \notag \\
&= d^{\times} ((x,e), (x,\gamma^{-1})) = d_{\Gamma} (e,\gamma^{-1}) = \vert \gamma^{-1} \vert = \vert \gamma \vert, \notag
\end{align}
where $d^{\times}$ stands for the max metric on the product $X \times \Gamma$.
\end{proof}

\begin{DefRef}{leftbdd44}
Let $b \colon X \to X^{bdd}$ be the natural map given by $b(x) = [x,e]$ in the orbit space $X \times_{\Gamma} \Gamma$.
\end{DefRef}

\begin{PropRef}{HJDSEO}
The map $b$ is a coarse map.
\end{PropRef}

\begin{proof}
Suppose $d^{bdd} ([x_1,e], [x_2,e]) \le D$, then $d^{\times} ((x_1,e), (\gamma x_2,\gamma)) \le D$ for some $\gamma \in \Gamma$,
so $d(x_1, \gamma x_2) \le D$ and $\vert \gamma \vert \le D$.
Since the left action of $\Gamma$ on $X^{bdd}$ is bounded, there is a function $B_{\ast}$ guaranteed by Lemma \refL{JUQASW2}.
Now
\[
d (x_1, x_2) \le d(x_1, \gamma x_2) + d (x_2, \gamma x_2)
\le D + B_{\ast} (D).
\]
This verifies that $b$ is proper.  It is clearly distance reducing, therefore uniformly expansive with $l(r)=r$.
\end{proof}

If we think of $X^{bdd}$ as the set $X$ with the metric induced from the bijection $b$, the map $b$ becomes the coarse identity map between the metric space $X$ with a left action of $\Gamma$ and the metric space $X^{bdd}$ where the \textit{action is made left-bounded}.

Of course we now have a construction we can apply in any instance of a free properly discontinuous action by $\Gamma$ on a metric space $X$.  In some cases $X$ has a very canonical set of metrics, for example in the case of a fundamental group acting cocompactly on the universal cover $X$.  On the other hand, if $\Gamma$ is the fundamental group of a manifold embedded in a Euclidean space, the normal bundle doesn't have such a set of metrics.  The authors are particularly interested in using left-bounded metrics on the normal bundles.  In these cases the map induced from $b$ preserves the most relevant $K$-theoretic information.

\bigskip

Finally, we describe one more coarse geometric construction that we will need.
For any metric space $X$, the new space $TX$ will be related to the cone
construction.

\begin{DefRef}{preTX}
Start with any set $Z$, let $S \subset Z \times Z$ denote any symmetric and
reflexive subset with the property that
\begin{itemize}
\item for any $z$, $z'$, there are elements $z_0$, $z_1$, ... , $z_n$ so that $z_0 = z$, $z_n = z_0$,
and $(z_i, z_{i+1}) \in S$.
\end{itemize}
Let $\rho \colon S \to \mathbb{R}$ be any function for which the following
properties hold:
\begin{itemize}
\item $\rho (z_1, z_2) = \rho (z_2, z_1)$ for all $(z_1, z_2) \in S$,
\item $\rho (z_1, z_2) = 0$ if and only if $z_1 = z_2$.
\end{itemize}

Given such $S$ and $\rho$, we may define a metric $d$ on $Z$ to be the largest metric $D$
for which $D (z_1, z_2) \le \rho (z_1, z_2)$ for all $(z_1, z_2) \in S$. This means that $d$ is given
by
\[
d(z_1, z_2) = \inf_{n,\{z_0,z_1,...,z_n\}}
\sum_{i=0}^n
\rho (z_i, z_{i+1}).
\]
\end{DefRef}

\begin{DefRef}{TX}
Let $k \ge 1$ be a real number or $\infty$.
We define a metric space $T_k X$ by first declaring that the underlying
set is $X \times \mathbb{R}$. Next, let $S$ to be the set consisting of pairs of
the form $((x, r), (x, r'))$ or of the form $((x, r), (x', r))$. We then define $\rho$ on $S$
by
\[
\rho ((x, r), (x, r')) = \vert r - r' \vert ,
\]
and
\[
\rho ((x, r), (x', r)) = \left\{
                           \begin{array}{ll}
                             d (x, x'), & \hbox{if $r \le 1$;} \\
                             r d (x, x'), & \hbox{if $1 \le r \le k$;} \\
                             k d (x, x'), & \hbox{if $k \le r$.}
                           \end{array}
                         \right.
\]
Since $\rho$ clearly satisfies the hypotheses of the above definition, we set the metric
on $T_k X$ to be $d_k$.  It is clear that the identity maps $T_k X \to T_m X$ are distance nondecreasing for $m \ge k$.

For simplicity of notation, we will write $TX$ for the colimit $T_{\infty} X$. 
We also extend the definition to pairs of metric spaces $(TX, TY)$,
where $Y$ is a subset of $X$ given the restriction of the metric on $X$.
\end{DefRef} 

\begin{PropRef}{TRT}
(1) If $X$ is a compact subset of a Euclidean space, one obtains
$\Knc (TX,R) = X_+ \wedge \Knc (R)$.

(2) $\Knc (T \mathbb{R}^n,R) \simeq \Sigma \Knc (\mathbb{R}^n,R) \simeq
\Sigma^{n+1} \Knc (R)$.
\end{PropRef}

\begin{proof}
	These are immediate consequences of applying standard bounded $K$-theory tools from \cite{mCeP:97,ePcW:89} to the cone construction. 
\end{proof}

\SecRef{Asymptotic transfer in parametrized \textit{K}-theory}{CPT}

Let $Y$ be a proper metric space.  Suppose $\Gamma$ is a finitely generated group which acts freely, properly discontinuously on $Y$.
In this case, one has the orbit space metric on $\Gamma \backslash Y$ as in Definition \refD{OrbitMet}.

Let $s \colon \Gamma \backslash Y \to Y$ be a section of the orbit space projection $p \colon Y \to \Gamma \backslash Y$.
We will require $s$ to be a coarse map.
For completeness, we provide an option for such map whenever $\Gamma \backslash Y$ is a finite simplicial complex.
This is the situation we will assume in the rest of the section.

\begin{LemRef}{HIMHIM}
Suppose $\Gamma \backslash Y$ is a finite simplicial complex $K$, then there is a section $s \colon \Gamma \backslash Y \to Y$ which is a coarse map.
\end{LemRef}

\begin{proof}
The complex $K$ is the union of finitely many open simplices which we denote $S_1, \ldots, S_t$.
Choose arbitrary base points $b_i$ in $S_i$.
Since each $S_i$ is contractible, an arbitrary map $s \colon \{ b_1, \ldots, b_t \} \to Y$ extends uniquely to a section $s \colon K \to Y$.
\end{proof}

Let $d$ be the standard simplicial metric in $Y$ which gives the simplicial orbit metric $d_{\Gamma \backslash Y}$ in $K$.

\begin{DefRef}{MKhhN22}
The image of $s$, $F \subset Y$, is a bounded fundamental domain for the action of $\Gamma$, since it is a union of finitely many simplices.

Let $\epsilon \ge 0$ and let $\Omega_{\epsilon}$ be the subset of $\Gamma$ given by
\[
\Omega_{\epsilon} = \{ \gamma \in \Gamma \mid d(F, \gamma F) \le \epsilon \}.
\]
Since the action of $\Gamma$ is properly discontinuous, the set $\Omega_{\epsilon}$ is finite.

Let $D_{\epsilon}$ be the maximal norm of an element in $\Omega_{\epsilon}$.
\end{DefRef}

\begin{PropRef}{POIU56}
If $d_{\Gamma \backslash Y} (z_1, z_2) \le \epsilon$ then there is an element $\omega$ of $\Omega_{\epsilon}$ such that $d(s(z_1), \omega s(z_2)) \le \epsilon$.
\end{PropRef}

\begin{proof}
Let $\omega$ be the element of $\Gamma$ such that
\[
d_{\Gamma \backslash Y} (z_1, z_2) = d(s(z_1), \Gamma s(z_2)) = d(s(z_1), \omega s(z_2)).
\]
Then clearly, if $d_{\Gamma \backslash Y} (z_1, z_2) \le \epsilon$ then $\omega \in \Omega_{\epsilon}$.
\end{proof}

Recall the construction of $T_k X$ from Definition \refD{TX}. 
A choice of $s$ determines a section $T_k s \colon T_k (\Gamma \backslash Y) \to T_k Y$.

\begin{PropRef}{POIU5}
$T_k s$ is a coarse map for all $k$.
\end{PropRef}

\begin{proof}
Explicitly, $(T_k s)(z,r) = (s(z),r)$.
So if $s$ is bounded by $D$ then $T_k s$ is bounded by $kD$ and so is coarse.
\end{proof}

We now consider the situation where $\Gamma$ acts on $(TY)^{bdd}$ via $\gamma (y,r) = (\gamma y, r)$, and $Tp \colon (TY)^{bdd} \to TK$ is the orbit space projection.
Let $D_1$ be the maximal norm of an element in $\Omega_1$.

\begin{PropRef}{POIU5p}
Suppose ${d}_{TK} ((z_1,r), (z_2,r)) \le R$ where ${d}_{TK}$ is the orbit metric in $TK$.
Then $d_{(TY)^{bdd}} ((s(z_1),r), (s(z_2),r)) \le R + 2D_1$.
\end{PropRef}

\begin{proof}
If $\gamma$ is an element of $\Gamma$ such that
\[
{d}_{TK} ((z_1,r), (z_2,r)) = d_{(TY)^{bdd}} ((s(z_1),r), (\gamma s(z_2),r)),
\]
then $\| \gamma \| = 1$
and
$d_{(TY)^{bdd}} ((\gamma s(z_2),r), (\gamma s(z_2),r)) \le 2 D_1$.  The result follows from the triangle inequality.
\end{proof}

Let $Ts \colon TK \to (TY)^{bdd}$ be given by $Ts (z,r) = (s(z), r)$.

\begin{CorRef}{POIU5p1}
$Ts$ is a coarse map and is a section of the projection $Tp$.
\end{CorRef}

\begin{proof}
Both facts follow from $Ts$ being the colimit of $T_k s$ and Proposition \refP{POIU5}.	
\end{proof}

\begin{NotRef}{GHYAS}
	We will denote the product metric space $\Gamma \times (TX)^{bdd}$ by $T_{\Gamma}(X)$.
\end{NotRef}

Notice that when the action of $\Gamma$ on $X$ is trivial then $(TX)^{bdd} = TX$, and so
in this case $T_{\Gamma}(X) = \Gamma \times TX$.

\begin{DefRef}{MKOERN1}
Given an isometric action of $\Gamma$ on $X$, there is the induced action on $\Gamma \times X \times \mathbb{R}$ given by $\gamma (\gamma',x,r) = (\gamma \gamma', \gamma x,r)$. This induces an isometric action on $T_{\Gamma}(X)$.
Recall that we are given a free, properly discontinuous isometric action of $\Gamma$ on a proper metric space $Y$.
Now we have the associated metric spaces $T_{\Gamma}(\Gamma \backslash Y)$ and $T_{\Gamma}(Y)$ with isometric actions by $\Gamma$.

Suppose $\Gamma \backslash Y$ is finite, then we can choose a coarse section $s \colon \Gamma \backslash Y \to Y$ of the projection $p \colon Y \to \Gamma \backslash Y$ as in Lemma \refL{HIMHIM}.
The sections $T s \colon T (\Gamma \backslash Y) \to T (Y)$ assemble to give a map
\[
T_{\Gamma} s \colon T_{\Gamma} (\Gamma \backslash Y) \longrightarrow T_{\Gamma} (Y)
\]
defined by
\[
Ts (\gamma, z, r) = (\gamma, Ts(z), r).
\]
\end{DefRef}

\begin{PropRef}{HIMHIOO}
The map $T_{\Gamma} s$ is a coarse map.
\end{PropRef}

\begin{proof}
Products of coarse maps are coarse.
\end{proof}

\begin{RemRef}{IOUY}
We note that if the action of $\Gamma$ on $T_{\Gamma} Y$ is via $\gamma' (\gamma, y, r) = (\gamma, \gamma'  y, r)$ then $T_{\Gamma} (\Gamma \backslash Y)$ is precisely the orbit space with the orbit space metric.
The projection $T_{\Gamma} p \colon T_{\Gamma} (Y) \to T_{\Gamma} (\Gamma \backslash Y)$ is the orbit space projection, and the map $Ts$ is a section of this projection.
\end{RemRef}

This is certainly different from the situation of main interest to us.  We are interested in the diagonal action of $\Gamma$ on $T_{\Gamma} Y$ given by $\gamma' (\gamma, y, r) = (\gamma' \gamma, \gamma' \, y, r)$.
We now construct and examine a section $s_{\Delta} \colon \Gamma \backslash T_{\Gamma} Y \to T_{\Gamma} Y$.

\begin{DefRef}{MKOERNN}
(1) As soon as a section $s \colon \Gamma \backslash Y \to Y$ is chosen, there is a well-defined function
$t \colon Y \to \Gamma$ determined by
$t(y) \, y = s([y])$.

(2) The diagonal action of $\Gamma$ on $T_{\Gamma} Y$ gives the orbit space projection
\[
p_{\Delta} \colon T_{\Gamma} (Y) \longrightarrow \Gamma \backslash T_{\Gamma} Y
\]
which endows $\Gamma \backslash T_{\Gamma} Y$ with the orbit space metric $d_{\Delta}$.
We define a map
\[
s_{\Delta} \colon \Gamma \backslash T_{\Gamma} Y \longrightarrow T_{\Gamma} (Y)
\]
by
$s_{\Delta} ([\gamma, y, r]) = (t(y) \, \gamma, s([y]), r)$.
Here $[y]$ stands for the class $p(y)$ and $[\gamma, y, r]$ for $p_{\Delta} (\gamma, y, r)$.
\end{DefRef}

\begin{NotRef}{KLOQA}
Given a subset $S \subset \Gamma$, let $T_S X$ denote the metric subspace $S \times X \times \mathbb{R}$ of $T_{\Gamma} X$.
\end{NotRef}

\begin{PropRef}{HIMHIOO2}
The map $s_{\Delta}$ is a section of $p_{\Delta}$.
It is not necessarily a coarse map.
However the restriction of $s_{\Delta}$ to each $T_S (\Gamma \backslash Y)$ for a bounded subset $S \subset \Gamma$ is coarse.
\end{PropRef}

\begin{proof}
By the defining property of $t(y)$,
\[
p_{\Delta} (\gamma, y, r) = p_{\Delta} (\gamma, t(y)^{-1} s([y]), r) = p_{\Delta} (t(y) \gamma, s([y]), r).
\]
So
$p_{\Delta} s_{\Delta} ([\gamma, y, r]) = p_{\Delta} (\gamma, y, r) = [\gamma, y, r]$.
A bound on the norm of $\gamma$ in $S$ gives a linear coefficient to exhibit the restriction of $s_{\Delta}$ to $T_S (\Gamma \backslash Y)$ as a coarse map.
In general, the bound for the action of $\gamma$ on $Y$ grows indefinitely.
\end{proof}

The function
\[
l \colon  T_{\Gamma} (\Gamma \backslash Y) \longrightarrow \Gamma \backslash T_{\Gamma} Y
\]
is given by $l (\gamma, z, r) = [\gamma, s(z), r]$.

\begin{PropRef}{HIMHPIMt}
$l$ is a coarse bijection.
\end{PropRef}

\begin{proof}
The function $l$ is proper and distance non-increasing, therefore coarse.
The inverse
$l^{-1} \colon \Gamma \backslash T_{\Gamma} Y \to T_{\Gamma} (\Gamma \backslash Y)$ is given by $l^{-1} ([\gamma, y, r]) = (t(y) \gamma, [y], r)$.
\end{proof}

\begin{DefRef}{poirs}
We define a function $u \colon T_{\Gamma} Y \to T_{\Gamma} Y $ by
\[
u(\gamma, y, r) = (t(y)^{-1} \gamma, y, r).
\]
\end{DefRef}

Clearly, $u$ is the identity on $\im (Ts)$ since $y$ in $\im (s)$ gives $t(y) = 1$ and $s([y]) = y$.

The following diagram shows the relationship between the geometric maps we have defined.
\[
\xymatrix{
 T_{\Gamma} Y   \ar@<1ex>[dd]^-{p_{\Delta}}
&&T_{\Gamma} Y  \ar@<1ex>[dd]^-{T_{\Gamma} p} \ar[ll]_-{u}\\
\\
 \Gamma \backslash T_{\Gamma} Y \ar@<1ex>[uu]^-{s_{\Delta}} \ar@<0.9ex>[rr]^-{l^{-1}}
&&T_{\Gamma} (\Gamma \backslash Y) \ar@<1ex>[uu]^-{T_{\Gamma} s} \ar@<0.5ex>[ll]^-{l}
}
\]

Let $B(k)$ stand for the subset of all elements $\gamma$ in $\Gamma$ with $\| \gamma \| \le k$.

\begin{DefRef}{poirs66}
Let $X$ and $Y$ be arbitrary metric spaces.
Given a function $f \colon T_{\Gamma} X \to T_{\Gamma} Y$, we say $f$ is $p$-\textit{bounded} if there is a function $\lambda \colon [0, \infty) \to [0, \infty)$ such that for all $k$ the restriction of $f$ to $T_{B(k)} X$ is a coarse map bounded by $\lambda (k)$.

It is clear that compositions of $p$-bounded maps are $p$-bounded.
\end{DefRef}

In general, this property is sufficient in order to induce an additive functor 
\[
\varrho \colon
\mathcal{B}^{\Gamma}_{p} ((TX)^{bdd})^{\Gamma} \longrightarrow
\mathcal{B}^{\Gamma}_{p} ((TY)^{bdd})^{\Gamma}.
\]
This can be seen by directly extending the proof of Theorem \refT{JARI}.
We will verify this in the case of the function
\[
s' \colon T_{\Gamma} (\Gamma \backslash Y) \to T_{\Gamma} Y
\]
defined as the composition $u \circ Ts$.

\begin{PropRef}{MNBTR}
The function $u$ is $p$-bounded.
The function $l^{-1}$ defined in the proof of Proposition \refP{HIMHPIMt} is $p$-bounded.
\end{PropRef}

\begin{proof}
Assuming $\| \gamma_1 \| \le k$, $\| \gamma_2 \| \le k$, and $d ((y_1,r_1), (y_2,r_2)) \le d$,
we estimate
\[
d((t(y_1)^{-1} \gamma_1, y_1, r_1), (t(y_2)^{-1} \gamma_2, y_2, r_2))
\le \min \{ d( t(y_1)^{-1} \gamma_1, t(y_2)^{-1} \gamma_2 ), d \}.
\]
Now
\begin{equation} \begin{split}
d(t(y_1)^{-1} \gamma_1, t(y_2)^{-1} \gamma_2) =\ & d(\gamma_1, t(y_1) t(y_2)^{-1} \gamma_2))\\
\le\ &d (\gamma_1 , \gamma_2) + d(\gamma_2, t(y_1) t(y_2)^{-1} \gamma_2))\\
\le\ &2k + \| t(y_1) t(y_2)^{-1} \|\\
\le\ &2k +d +2D_1,\\
\end{split} \notag \end{equation}
where $D_1$ is the constant from Proposition \refP{POIU5p}.
The second statement is proved similarly.
\end{proof}

Since $s' \colon T_{\Gamma} (\Gamma \backslash Y) \to T_{\Gamma} Y$ is the composition $u \circ Ts$, it is therefore $p$-bounded.

\begin{PropRef}{MNBTrrr}
We have $s_{\Delta} = u \circ Ts \circ l^{-1}$.
\end{PropRef}

\begin{proof}
Since $l^{-1} ([\gamma, y, r]) = (t(y) \gamma, [y], r)$, we have
\[
T(s) \, l^{-1} ([\gamma, y, r]) = (t(y) \gamma, s([y]), r).
\]
Then also $u \, T(s) \, l^{-1} ([\gamma, y, r]) = (t(y) \gamma, s([y]), r)$ because $t(s[y])^{-1} = 1$.
\end{proof}

This proves that the diagram shown above is commutative.
It also follows that $s_{\Delta} = s \circ l^{-1}$ is $p$-bounded.

We are ready to define the parametrized transfer.  Suppose $\Gamma$ acts on $Y$ by isometries.

\begin{DefRef}{HIMHIMt}
The \textit{parametrized transfer map}
\[
P \colon K^{\Gamma}_{i} (T (\Gamma \backslash Y))^{\Gamma} \longrightarrow
K^{\Gamma}_{p} ((TY)^{bdd})^{\Gamma}
\]
is obtained from the $p$-bounded section $s \colon T_{\Gamma} (\Gamma \backslash Y) \to T_{\Gamma} (Y)$ via the induced additive functor
\[
\varrho \colon
\mathcal{B}^{\Gamma}_{i} (T (\Gamma \backslash Y))^{\Gamma} \longrightarrow
\mathcal{B}^{\Gamma}_{p} ((TY)^{bdd})^{\Gamma}.
\]
Here the action of $\Gamma$ on the set $T_{\Gamma} (\Gamma \backslash Y) = \Gamma \times \Gamma \backslash Y \times \mathbb{R}$ is given by $\gamma (\gamma',z,r) = (\gamma \gamma', z, r)$, and the action on
$T_{\Gamma} Y = \Gamma \times Y \times \mathbb{R}$ by $\gamma (\gamma',y, r) = (\gamma \gamma', \gamma y, r)$.

In order to describe how the map $s$ induces the functor $\varrho$, let us reiterate the interpretation in Definition \refD{EqCatDef} of an object of $\mathcal{B}^{\Gamma}_{i} (T (\Gamma \backslash Y))^{\Gamma}$.
It is a pair
$(F,\psi)$ where $F$ is an object of $\mathcal{B} (T_{\Gamma} (\Gamma \backslash Y))$ and $\psi$ is a function on
$\Gamma$ with $\psi_{\gamma} \in \Hom (F, \gamma F)$ in $\mathcal{B} (T_{\Gamma} (\Gamma \backslash Y))$ such that $\psi_{\gamma}$ is of filtration $0$ when projected to $\Gamma$,
$\psi_e = \id$, and $\psi_{\gamma_1 \gamma_2} =
\gamma_1 \psi_{\gamma_2}  \psi_{\gamma_1}$.

Given $(F,\psi)$ in $\mathcal{B}^{\Gamma}_{i} (T (\Gamma \backslash Y))^{\Gamma}$, we proceed to define $(\varrho F, \varrho \psi)$ in $\mathcal{B}^{\Gamma}_{p} ((TY)^{bdd})^{\Gamma}$.

Using the notation $[y]$ for the orbit of $y \in Y$, we set $(\varrho F)_{(\gamma,y,r)} = F_{(\gamma,[y],r)}$ if $y=\gamma s[y]$ and $0$ otherwise. 
Concretely, we tensor a lift of the module $F_{(\gamma,[y],r)}$ with the action of $\Gamma$.
The morphisms $(\varrho \psi)_{\gamma} \in \Hom (\varrho F, \gamma \varrho F)$ are defined
by making the component
\[
((\varrho \psi)_{\gamma})_{(\gamma_1,y_1,r_1), (\gamma_2,y_2,r_2)} \colon (\varrho F)_{(\gamma_1,y_1,r_1)} \longrightarrow (\gamma \varrho F)_{(\gamma_2,y_2,r_2)},
\]
where $(\varrho F)_{(\gamma_1,y_1,r_1)} = F_{(\gamma_1,[y_1],r_1)}$ and
\[
(\gamma \varrho F)_{(\gamma_2,y_2,r_2)} = (\varrho F)_{(\gamma^{-1} \gamma_2, \gamma^{-1} y_2, r_2)} = F_{(\gamma^{-1} \gamma_2, [\gamma^{-1} y_2], r_2)} = F_{(\gamma^{-1} \gamma_2, [y_2], r_2)},
\]
the $0$ homomorphism unless
\[
t(y_1) \, t(y_2)^{-1} = \gamma_1^{-1} \gamma_2.
\]
In the latter case, the component $((\varrho \psi)_{\gamma})_{(\gamma_1,y_1,r_1), (\gamma_2,y_2,r_2)}$ is identified with
\[
(\psi_{\gamma})_{(\gamma_1, [y_1],r_1), (\gamma_2, [y_2],r_2)} \colon F_{(\gamma_1,[y_1],r_1)} \longrightarrow (\gamma F)_{(\gamma_2, [y_2],r_2)} = F_{(\gamma^{-1} \gamma_2, [y_2], r_2)}.
\]
The point is that in this case there is an element $\gamma' \in \Gamma$ such that $t(y_1) = \gamma_1^{-1} (\gamma')^{-1}$ and $t(y_2) = \gamma_2^{-1} (\gamma')^{-1}$, so $y_1 = \gamma' \gamma_1 s([y_1])$ and $y_2 = \gamma' \gamma_2 s([y_2])$.

This shows that  if $\psi_{\gamma}$ is a morphism of $\mathcal{B} (T_{\Gamma} (\Gamma \backslash Y))$ then $(\varrho \psi)_{\gamma}$ is a morphism of $\mathcal{B}^{\Gamma}_{p} ((TY)^{bdd})^{\Gamma}$.
Also $(\varrho \psi)_{e} = \psi_e = \id$, and clearly $(\varrho \psi)_{\gamma_1 \gamma_2} =
\gamma_1 (\varrho \psi)_{\gamma_2}  (\varrho \psi)_{\gamma_1}$.

If a morphism $\phi \colon (F, \psi) \to (F', \psi')$ in $\mathcal{B}^{\Gamma}_i (T(\Gamma \backslash Y))^{\Gamma}$ is given by $\phi \colon F \to F'$ in $\mathcal{B} (T_{\Gamma} (\Gamma \backslash Y))$, then $\varrho \phi \colon \varrho F \to \varrho F'$ is defined by
\[
(\varrho \phi)_{(\gamma_1, y_1,r_1), (\gamma_2, y_2,r_2)} =
\begin{cases}
0 &\text{if $t(y_1) \, t(y_2)^{-1} \ne \gamma_1^{-1} \gamma_2$,} \\
\phi_{(\gamma_1, [y_1], r_1), (\gamma_2, [y_2], r_2)} &\text{otherwise}.
\end{cases}
\]
It is a morphism of $\mathcal{B}_{p} ( (TY)^{bdd})$ because the section $s$ is $p$-bounded, and the action of $\Gamma$ is bounded.
\end{DefRef}

One useful idea when computing the parametrized transfer map is a decomposition of the map as a homotopy colimit of the transfers associated to subspaces of $\Gamma \backslash Y$.  Whether the transfer is equivalent to the homotopy colimit depends on the circumstances, but the  colimit itself can be always constructed as follows.

Suppose $\Gamma \backslash Y$ is the union of subspaces $S_1$ and $S_2$.  We will denote the covering of $\Gamma \backslash Y$ by $S_1$, $S_2$, and $S_1 \cap S_2$ by $\mathcal{S}$.   We have seen there are spectra associated to full subcategories $\mathcal{B}^{\Gamma}_i (T(\Gamma \backslash Y))^{\Gamma}_{<S}$ of $\mathcal{B}^{\Gamma}_i (T(\Gamma \backslash Y))^{\Gamma}$ for $S$ in $\mathcal{S}$.  We also know that in this situation, when the action on $Y$ is trivial, there are equivalences $\mathcal{B}^{\Gamma}_i (T(\Gamma \backslash Y))^{\Gamma}_{<S} \simeq \mathcal{B}^{\Gamma}_i (TS)$. Now there is $K^{\Gamma}_{i} (T \mathcal{S})^{\Gamma}$ defined as the homotopy colimit of the diagram 
\[
\xymatrix{
 K^{\Gamma}_{i} (T (S_1 \cap S_2))^{\Gamma} \ar[r] \ar[d]
&K^{\Gamma}_{i} (T S_1)^{\Gamma} \\
 K^{\Gamma}_{i} (T S_2)^{\Gamma}
&
}
\]
and we know from Theorem \refT{Exc2} that, if $S_1$ and $S_2$ are a coarsely antithetic pair, then $K^{\Gamma}_{i} (T \mathcal{S})^{\Gamma}$ is equivalent to $K^{\Gamma}_{i} (T (\Gamma \backslash Y))^{\Gamma}$.

If we denote the full preimage of $S$ in $Y$ along the projection $p$ by $\widetilde{S}$, there is a parametrized transfer map 
\[
P_S \colon K^{\Gamma}_{i} (T S)^{\Gamma} \longrightarrow
K^{\Gamma}_{p} ((T\widetilde{S})^{bdd})^{\Gamma}
\]
for each subspace $S$ from the covering $\mathcal{S}$.  

It is easy to see that if $\mathcal{S}$ is coarsely antithetic, the covering of $TY$ by $T\widetilde{S}$ is again coarsely antithetic.  In this case, the colimit $K^{\Gamma}_{p} ((T\widetilde{\mathcal{S}})^{bdd})^{\Gamma}$ of 
\[
\xymatrix{
 K^{\Gamma}_{p} ((T \widetilde{S}_1 \cap T \widetilde{S}_2)^{bdd})^{\Gamma} \ar[r] \ar[d]
&K^{\Gamma}_{p} ((T \widetilde{S}_1)^{bdd})^{\Gamma} \\
 K^{\Gamma}_{p} ((T \widetilde{S}_2)^{bdd})^{\Gamma}
&
}
\]
is equivalent to $K^{\Gamma}_{p} ((TY)^{bdd})^{\Gamma}$ by Theorem \refT{BDDEXCI25}. Using the notation from the last two paragraphs and induction, we obtain the following fact.

\begin{ThmRef}{BVCZ}
Suppose $\mathcal{S}$ is a finite covering of $\Gamma \backslash Y$.  Then there is a parametrized transfer map 
\[
P_{\mathcal{S}} \colon 
K^{\Gamma}_{i} (T \mathcal{S})^{\Gamma}
\longrightarrow
K^{\Gamma}_{p} ((T\widetilde{\mathcal{S}})^{bdd})^{\Gamma}
\]
which is the colimit of the maps $P_S$.  If $\mathcal{S}$ is a coarsely antithetic covering, in the sense that all pairs of subsets are mutually coarsely antithetic, then $P_{\mathcal{S}}$ is the 
	parametrized transfer map $P_{\Gamma \backslash Y}$.
\end{ThmRef}

The following is a geometric situation where the transfer is particularly useful.  Let ${N}$ be the total space of the normal disk bundle.
Then the discussion in this section applies to $Y = \widetilde{N}$ and $\Gamma \backslash Y = {N}$.

\begin{ThmRef}{UICSR}
There exists a map of nonconnective spectra
\[
P \colon K^{\Gamma}_i (TN)^{\Gamma} \longrightarrow K^{\Gamma}_{p} ( (T\widetilde{N})^{bdd})^{\Gamma}
\]
and its relative version 
\[
P' \colon K^{\Gamma}_i (T {N}, T \partial {N})^{\Gamma} \longrightarrow K^{\Gamma}_{p} ((T \widetilde{N})^{bdd}, (T \partial \widetilde{N})^{bdd})^{\Gamma}.
\]
Both are right inverses to the maps induced by the covering projection $p \colon \widetilde{N} \to N$.
\end{ThmRef}

This theorem can be viewed as an analogue to the asymptotic transfers and their properties that appear, for example, in section 12 of \cite{aBtFlJhR:04} or Theorem 6.1 of \cite{aBwLhR:08a}. 

\begin{proof}
We obtain $P$ as a parametrized transfer map from Definition \refD{HIMHIMt}.  The functor $\varrho$ induces an additive functor 
\[
\varrho' \colon
\mathcal{B}^{\Gamma}_i (T {N}, T \partial {N})^{\Gamma}
\longrightarrow
\mathcal{B}^{\Gamma}_{p} ((T \widetilde{N})^{bdd}, (T \partial \widetilde{N})^{bdd})^{\Gamma}
\]
between Karoubi quotients as in Proposition \refP{FibPair}.
This gives the relative version $P'$.

Recall that $\varrho$ was defined by the assignment $(\varrho F)_{(\gamma,y,r)} = F_{(\gamma,[y],r)}$ if $y=\gamma s[y]$, so $P$ is a right inverse to the map induced from the functor $\varpi$ defined on objects by $(\varpi F)_{(\gamma,[y],r)} = F_{(\gamma,\gamma s[y],r)}$.   The latter is precisely the functor induced by $p$ since $[y]=p( \gamma s[y])$.
\end{proof}


\begin{thebibliography}{99}

\bibitem{sA:16}
S. Arnt, \textit{Fibred coarse embeddability of box spaces and proper isometric affine actions on $L^p$ spaces}, Bull. Belg. Math. Soc. Simon Stevin \textbf{23} (2016), 21--32.

\bibitem{aBtFlJhR:04}
A. Bartels, T. Farrell, L. Jones, and H. Reich, 
\textit{On the isomorphism conjecture in algebraic K–theory}, Topology \textbf{43} (2004), 157--213.

\bibitem{aBwLhR:08a}
A.~Bartels, W.~L\"{u}ck, and H.~Reich, 
\textit{The K-theoretic Farrell-Jones Conjecture for hyperbolic groups},
Invent. Math. (2008), 29--70.

\bibitem{aBwL:12}
A.~Bartels and W.~L\"{u}ck
\textit{The Borel Conjecture for hyperbolic and CAT(0)-groups},
Ann. Math. 175 (2012), 631--689.

\bibitem{jBdG:99}
J.C. Becker and D.H. Gottlieb, 
\textit{A history of duality in algebraic topology},
in History of topology, North-Holland, Amsterdam, 1999, 725-745.

\bibitem{mCeP:97}
{M. Cardenas and E.K. Pedersen},
{\it On the Karoubi filtration of a category},
$K$-theory, {\bf 12} (1997), 165--191.

\bibitem{gC:93}
{G. Carlsson}, 
{\it Proper homotopy theory and transfers for infinite groups}, 
in Algebraic Topology and its Applications,
MSRI Publications \textbf{27} Springer--Verlag (1993), 1--10.

\bibitem{gC:95}
\bysame,
{\it Bounded $K$-theory and the assembly map
in algebraic $K$-theory},
in {\it Novikov conjectures, index theory and rigidity},
{\it Vol. 2}
(S.C. Ferry, A. Ranicki, and J. Rosenberg, eds.),
Cambridge U. Press (1995), 5--127.

\bibitem{gCbG:04}
{G.~Carlsson and B.~Goldfarb}, \textit{The integral K-theoretic Novikov conjecture for
groups with finite asymptotic dimension}, Inventiones Math. {\bf
157} (2004), 405--418.

\bibitem{gCbG:00}
\bysame,
\textit{Controlled algebraic $G$-theory, I}, J. Homotopy Relat. Struct. \textbf{6} (2011), 119--159.

\bibitem{gCbG:15}
\bysame, \textit{On modules over infinite group rings}, Int. J. Algebra Comput. \textbf{26} (2016), 1--16.

\bibitem{gCbG:19}
\bysame, \textit{Bounded $G$-theory with fibred control},  
J. Pure Appl. Algebra \textbf{223} (2019), 5360--5395.

\bibitem{gCbG:20}
\bysame, \textit{Excision in equivariant fibred \textit{G}-theory}. \texttt{arXiv:1911.09397}

\bibitem{gCeP:93}
{G. Carlsson and E.K. Pedersen}, {\it Controlled algebra and the
Novikov conjecture for K- and L-theory}, Topology {\bf 34}
(1993), 731--758.

\bibitem{xCqWxW:13}
X. Chen, Q. Wang, and X. Wang, \textit{Characterization of the Haagerup property by fibred coarse embedding into Hilbert space}, Bull. London Math. Soc. \textbf{45} (2013), 1091--1099.

\bibitem{xCqWgY:13}
X. Chen, Q. Wang, and G. Yu,
\textit{The maximal coarse Baum-Connes conjecture for spaces which admit a fibred coarse embedding into Hilbert space}, Adv. Math. \textbf{249} (2013), 88--130.

\bibitem{fFlJ:86}
F.T. Farrell and L.E. Jones,
\textit{K-theory and dynamics. I}, 
Ann. Math. \textbf{124} (1986), 531--569.

\bibitem{mF:14}
M. Finn-Sell, 
\textit{Fibred coarse embeddings, a-T-menability and the coarse analogue of the Novikov conjecture}, 
J. Funct. Anal. \textbf{267} (2014), 3758--3782.

\bibitem{pMjS:06}
J.P. May and J. Sigurdsson,
\textit{Parametrized homotopy theory}, Amer. Math. Soc. (2006).

\bibitem{mMhS:18}
M. Mimura and H. Sako, 
\textit{Group approximation in Cayley topology and coarse geometry, Part II: Fibered coarse embeddings},
preprint, 2018.
\texttt{arXiv:1804.10614}

\bibitem{eP:84}
{E.K. Pedersen},
{\it On the $K_{-i}$-functors},
J. Algebra \textbf{90} (1984), 461--475.

\bibitem{ePcW:85}
{E.K. Pedersen and C. Weibel},
{\it A nonconnective delooping of algebraic $K$-theory},
in {\it Algebraic and geometric topology}
(A. Ranicki, N. Levitt, and F. Quinn, eds.),
Lecture Notes in Mathematics {\bf 1126},
Springer-Verlag (1985), 166--181.

\bibitem{ePcW:89}
{\bysame}, {\it $K$-theory homology of spaces}, in {\it Algebraic
topology} (G.~Carlsson, R.L.~Cohen, H.R.~Miller, and
D.C.~Ravenel, eds.), Lecture Notes in Mathematics {\bf 1370},
Springer-Verlag (1989), 346--361.

\bibitem{jR:03}
{J. Roe},
{\it Lectures on coarse geometry}, University Lecture Series, vol.~31, American Mathematical Society, 2003.

\bibitem{rT:83}
{R.W. Thomason},
{\it The homotopy limit problem},
in {\it Proceedings of the Northwestern homotopy theory conference}
(H.R. Miller and S.B. Priddy, eds.),
Contemp. Math. {\bf 19} (1983), 407--420.

\bibitem{kW:10}
{K.~Whyte}, \textit{Coarse bundles}, preprint, 2010. \texttt{arXiv:1006.3347}

\bibitem{gY:98}
{G.~Yu}, \textit{The Novikov conjecture for groups with finite asymptotic dimension}, Annals of Math. \textbf{147} (1998), 325--355.

\bibitem{gY:00}
\bysame, \textit{The coarse Baum-Connes conjecture for spaces which admit a uniform embedding into Hilbert space}, Invent. Math. \textbf{139} (2000), 201--240 
\end{thebibliography}
\end{document}